\documentclass[review]{elsarticle}

\usepackage{lineno,hyperref}

\usepackage{amsmath}
\usepackage{amssymb}
\usepackage{booktabs}
\usepackage{caption}
\usepackage{graphicx}
\usepackage[title]{appendix}
\usepackage{float}
\usepackage{subcaption}

\newtheorem{remark}{Remark}[section]
\newtheorem{thm}{Theorem}[section]
\newtheorem{lemma}{Lemma}[section]
\newtheorem{hypothesis}{Hypothesis}[section]
\newtheorem{corollary}{Corollary}[section]
\newtheorem{property}{Property}[section]
\newtheorem{definition}{Definition}[section]
\newtheorem{proposition}{Proposition}[section]
\numberwithin{equation}{section}

\begin{document}

\begin{frontmatter}
\hypersetup{linkcolor = blue,anchorcolor =red,citecolor = blue,filecolor = red,urlcolor = red,
            pdfauthor=author}


\title{Positive effects of multiplicative noise on the explosion of nonlinear fractional stochastic differential equations  \tnoteref{t1,t2}}
\tnotetext[t1]{\text{AMS Subject Classification}: 91D30 \sep 05C82 \sep 05C80 \sep 05C85.}
\tnotetext[t2]{This work is supported by the State Key Program of National Natural Science of China under Grant No.91324201. This work is also supported by the Fundamental Research Funds for the Central Universities of China under Grant 2018IB017, Equipment Pre-Research Ministry of Education Joint Fund Grant 6141A02033703 and the Natural Science Foundation of Hubei Province of China under Grant 2014CFB865.
}


\author[mymainaddress]{Fei Gao\corref{mycorrespondingauthor}}
\ead{gaof@whut.edu.cn}
\author[mainaddress]{Xinyi Xie}
\ead{2503551142@qq.com}
\author[mymainaddress]{Hui Zhan}
\ead{2432593867@qq.com}

\cortext[mycorrespondingauthor]{Corresponding author.}

\address[mymainaddress]{Department of Mathematics and Center for Mathematical Sciences, Wuhan University of Technology, Wuhan, 430070, China}
\address[mainaddress]{Department of Statistics and Center for Mathematical Sciences, Wuhan University of Technology, Wuhan, 430070, China}

\begin{abstract}
For the nonlinear stochastic partial differential equation which is driven by multiplicative noise of the form
\[D_t^\beta u = \left[ { - {{\left( { - \Delta } \right)}^s}u + \zeta \left( u \right)} \right]dt + A\sum\limits_{m \in Z_0^d} {\sum\limits_{j = 1}^{d - 1} {{\theta _m}{\sigma _{m,j}}\left( x \right)} }  \circ dW_t^{m,j},\;\; s \ge 1,\;\;\frac{1}{2} < \beta  < 1,\]
where $D_{t}^\beta $ denotes the Caputo derivative, $A>0$ is a constant depending on the noise intensity, $\circ$ represent the Stratonovich-type stochastic differential, we consider the blow-up time of its solutions. We find that the introduction of noise can effectively delay the blow-up time of the solution to the deterministic differential equation when $\zeta$ in the above equation satisfies some assumptions. A key element in our construction is using the Galerkin approximation and a priori estimates methods to prove the existence and uniqueness of the solutions to the above stochastic equations, which can be regarded as the fractional order extension of the conclusions in \cite{flandoli2021delayed}. We also verify the validation of hypotheses in the time fractional Keller-Segel and time fractional Fisher-KPP equations in 3D case.

\end{abstract}


\end{frontmatter}


\section{Introduction}
We consider the blow-up problem of the solutions of the fractional differential equation
\begin{equation}\label{1.1}
D_t^\beta u = \left[ { - {{\left( { - \Delta } \right)}^s}u + \zeta \left( u \right)} \right]dt,\quad \;{\rm{in}}\;(0,T)\times{\mathbb T^d} = {{{\mathbb R^d}} \mathord{\left/
 {\vphantom {{{\mathbb R^d}} {{\mathbb Z^d}}}} \right.
 \kern-\nulldelimiterspace} {{\mathbb Z^d}}}
\end{equation}
where $0<\beta<1$, $s \ge 1$, $d \ge 2$, and $T>0$,  $\Delta$ is the periodic Laplacian operator, with initial data $u\left( {0,x} \right) = {u_0}$. The precise definition of the Caputo derivative is given in Appendix \ref{appendixA}. We are mainly interested in fractional order equations, since
these models can preserve the genetic and memory properties of functions in practical problems, and the physical meaning of parameters also more explicit, cf.\cite{RN89,RN90,RN85,tien2013fractional,henry2002existence,glockle1995fractional}. For $\beta=1$ we recover the classical  nonlinear time fractional diffusion equations, whose theory is well known, cf. \cite{escudero2006fractional}.

However, for different initial data the  properties of the solutions to such equations is different, that is, the solutions may exist globally or blow up in finite time. The phenomenon that the solutions of ordinary or partial differential equations diverge in finite time is known as the blow up of solutions \cite{hirota2006numerical}, which exists in various fields, such as chemotaxis in biology, curvature flow in geometry, and fluid mechanics, see \cite{biler2010blowup,angenent1995asymptotic,fu2009well}. This problem has recently attracted some attention in mathematical analysis, such as the sufficient blow-up conditions for semi-discrete partial differential equations \cite{abia1996blow,abia1998blow}, the upper bound for the blow-up time of solutions of fully discrete partial differential equations \cite{abia2001euler}, the relationship between the blow-up set of semi-discrete equations and the continuous equations \cite{groisman2001asymptotic}.

We refer to \cite{zhang2015blow,RN81} for the basic theory of the blow-up problem for solutions of deterministic equations. In particular, the equation is called the fractional Fisher-KPP equation for the nonlinear part is taken as $ - u\left( {1 - u} \right)$. Alsaedi \cite{alsaedi2021global} also describes the asymptotic behavior of its bounded solutions. Recently, the relationship between mild and weak solutions has been considered in \cite{li2021blow}. Related literature is also mentioned in these papers.

In practical problems, complex systems often have a large number of noise perturbations, and it is well known that the properties of the solution may be affected by some perturbations. The resulting stochastic models have attracted the interest of many academics. For example, the suitably non-degenerate additive noise can improve the pathwise non-uniqueness for some deterministic equation\cite{da2010pathwise}. The strong uniqueness of a class of SPDE with a non-Lipschitz drift and an additive noise was studied in \cite{butkovsky2019regularization}, which is also called the regularization by noise. Multiplicative noise makes the solution of the stochastic transport linear equation converge to the deterministic equation when appropriate assumptions are satisfied \cite{galeati2020convergence}, and noise can improve the well-posedness for stochastic scalar conservation laws \cite{gess2018well}.

The main purpose of this paper is to discuss whether the presence of noise has an effect on the blow-up time of the solution, so we add multiplicative noise to the deterministic equation \eqref{1.1}, shaped as
\begin{equation} \label{GrindEQ__1_1_} 
D_t^\beta u = \left[ { - {{\left( { - \Delta } \right)}^s}u + \zeta \left( u \right)} \right]dt + A\sum\limits_{m \in Z_0^d} {\sum\limits_{j = 1}^{d - 1} {{\theta _m}{\sigma _{m,j}}\left( x \right)} }  \circ dW_t^{m,j}  \;\;\;\;\;\; {\rm{in}}\;\left( {0,T} \right) \times {\mathbb T^d},
\end{equation} 
where $\frac{1}{2} < \beta  < 1$, ${\mathbb T^d} = {\mathbb R^d}/{\mathbb Z^d}$ and $T>0$, $\circ$ represent the Stratonovich-type stochastic differential. We are mainly interested in the case where the function $\zeta$ is nonlinear. For $\beta =1$ we recover the stochastic partial differential equation (that will be shortened as SPDE), whose theory is widely studied, cf.\cite{flandoli2021delayed}. We also assume that we are given initial data
\begin{equation} \label{GrindEQ11}
u\left( {0,x} \right) = {u_0}
\end{equation}
where in principle ${u_0} \in {L^2}\left( {{\mathbb T^d}} \right)$. for reader’s convenience we have listed in Appendix \ref{appendixB} the most relevant sources to the applications that we know of.

We remind that Galerkin approximations and a prior estimates are efficient tools for verifying the uniqueness and existence of solutions, see in that respect the work of Blomker \cite{blomker2013galerkin} and Alikhanov \cite{alikhanov2010priori}. Therefore we use these methods in the following to prove the existence and uniqueness of the solution of the equation \eqref{GrindEQ__1_1_}-\eqref{GrindEQ11}.

\noindent\textbf{Outline of the paper and main results} First, some preliminary divisions of the blow-up time with respect to the solution of the equation \cite{galaktionov2002problem}. The case is called no blow-up: exists global solutions, which remain consistent and bounded in time. The second case is that the global solution explodes in infinite time. The standard blow-up case we know it is that the solution blows up in a finite time. When the explosion occurs at the moment $t = 0$, it is known as instantaneous blow-up. Additionally, the type of equation in which all solutions explode in finite time is also known as the Fujita problem \cite{RN91}. In this paper we mainly consider the standard blow-up case.

In Section \ref{sec:3}, we first introduce the meaning of each parameter in the model \eqref{GrindEQ__1_1_}, and make assumptions on the nonlinear part. Then two important conclusions of this paper are given. Theorem \ref{theorem3.1} shows that the solution of the system \eqref{GrindEQ__1_1_}-\eqref{GrindEQ11} exists in finite time when certain conditions are satisfied, namely, it does not blow up. Theorem \ref{theorem3.2}, on the other hand, states that the perturbation of multiplicative noise allows the existence of the solution to be extended from finite time to infinity, namely, the solution explosion time can be effectively delayed.

Section \ref{sec:4} is devoted to quote two examples, the fractional Keller-Segel and fractional Fisher-Kpp models, in which we verify the validity of the assumptions on nonlinear terms in Section \ref{sec:3}. Moreover, we discuss the initial value conditions that make the solution explode in finite time in the fractional Fisher-Kpp system.

Section \ref{sec:5} studies the global solution of equation \eqref{GrindEQ__1_1_}-\eqref{GrindEQ11} by introducing the cut-off function. We first discuss the boundedness of the solution in finite dimensions, using galerkin approximations and a prior estimate to prove the existence of uniqueness of the solution. And finally, we give the complete proof of Theorem \ref{theorem3.1} and Theorem \ref{theorem3.2} in Section \ref{proof}.

In Appendix A, we have collected the definitions of fractional order  derivative, as well as some definitions and lemmas that are significant in the process of proof.

\vspace{1cm}
\noindent{\scshape Notations.} The norm of $f \in {L^p}\left( {0,T;{H^s}\left( {{\mathbb T^d}} \right)} \right)$ is abbreviated as ${\left\| f \right\|_{{L^p}{H^s}}}$, the symbols ${\left\| f \right\|_{{L^\infty }{L^2}}},{\left\| f \right\|_{{C^\kappa }{H^s}}}$ and others have the same meaning. For any $f,g \in {L^2}\left( {{\mathbb T^d}} \right)$, the inner product is denoted as $\left\langle {f,g} \right\rangle$. For the prevention of symbol abuse, we also have the symbol $\left\langle {f,g} \right\rangle$ for the dual product of $f\in H^{s} \left({\rm {\mathbb T}}^{d} \right)$ and $g\in H^{-s} \left({\rm {\mathbb T}}^{d} \right)$. The notation $\left[ {\; \cdot \;,\; \cdot \;} \right]$ means the quadratic covariation process. We shall call $s -$Laplacian of $u$ the function $- {\left( { - \Delta } \right)^s}u$, this is consistent with the use in the standard case $s=1$.

In the following representations, the symbol “ $\lesssim$ ” indicates that the number and equation on the left is less than or equal to $K$ times the number and equation on the right, where $K$ is a constant greater than 0, and its value will change under different circumstances. The dependence of the constants on parameters will be clearly written only when necessary.

\section{Models and main results}
\label{sec:3}

We will introduce the meaning of each parameter in model \eqref{GrindEQ__3_50_},
\begin{equation}\label{GrindEQ__3_50_}
D_t^\beta u = \left[ { - {{\left( { - \Delta } \right)}^s}u + \zeta \left( u \right)} \right]dt,
\end{equation}
note that $\Delta $ represents the periodic Laplacian operator, $s\ge 1$ is fixed, $\frac{1}{2} <\beta <1$, and $\zeta \left( u \right)$ is a nonlinear function with respect to $u$. The above equation is defined on the tours ${\mathbb T^d} = {{{\mathbb R^d}} \mathord{\left/{\vphantom {{{\mathbb R^d}} {{\mathbb Z^d}}}} \right.\kern-\nulldelimiterspace} {{\mathbb Z^d}}}$, in this paper we discuss the high dimensional case, where $d \ge 2$. In order to make the conclusions more general, we assume that the nonlinear function $\zeta$ satisfies some assumptions that are motivated by \cite{liu2015stochastic} and \cite{flandoli2021delayed}. We will test the validity of this hypothesis with some classical equations in Section \ref{sec:4}.

\newcommand{\upcite}[1]{\textsuperscript{\textsuperscript{\cite{#1}}}}

\begin{hypothesis}\label{hypothesis3.1} \upcite{flandoli2021delayed} The unknown nonlinear term $\zeta$ satisfies the assumptions (\romannumeral1)-(\romannumeral4).

(\romannumeral1) The unknown nonlinear function $\zeta$ is a continuous mapping from Sobolev space $H^{s-\gamma _{1} }$ to Sobolev space $H^{-s} $, and the Sobolev norm of the function $\zeta$ satisfies the inequality
\[{\left\| {\zeta \left( u \right)} \right\|_{{H^{ - s}}}} \le C\left( {1 + \left\| u \right\|_{{L^2}}^{{a_1}}} \right)\left( {1 + {{\left\| u \right\|}_{{H^s}}}} \right),\]
\noindent with parameters $a_{1} \ge 0$ and $0 < {\gamma _1} < s$.

(\romannumeral2) The inner product satisfies the inequality
\[\left|\left\langle \zeta\left(u\right),u\right\rangle \right|\le C \left(1+\left\| u\right\| _{L^{2} }^{a_{2} } \right)\left(1+\left\| u\right\| _{H^{s} }^{\gamma _{2} } \right),\] 
\noindent with parameters $a_{2} \ge 0$, $0 < {\gamma _2} < 2$. More generally, scaling ${L^2}$-norm of $u\left( t \right)$ further by the interpolation inequality gives
\[\begin{array}{l}
\left| {\left\langle {\zeta \left( u \right),u} \right\rangle } \right| \le C\left( {1 + \left\| u \right\|_{{H^s}}^{{\gamma _2}}} \right)\left( {1 + \left\| u \right\|_{{H^s}}^{{a_2}\frac{\delta }{{s + \delta }}}} \right)\left( {1 + \left\| u \right\|_{{H^{ - \delta }}}^{{a_2}\frac{s}{{s + \delta }}}} \right)\\
{\rm{\;\;\;\;\;\;\;\;\;\;\;\;\;\;\;\;\;            }} \le C \left( {1 + \left\| u \right\|_{{H^s}}^{{\gamma _2} + {a_2}\frac{\delta }{{s + \delta }}}} \right)\left( {1 + \left\| u \right\|_{{H^{ - \delta }}}^{{a_2}\frac{s}{{s + \delta }}}} \right),
\end{array}\]
where $\delta$ can be taken to be small enough.

(\romannumeral3) Extending the first inequalities in condition (\romannumeral2), it holds
\[\left|\left\langle u-v,\zeta\left(u\right)-\zeta\left(v\right)\right\rangle \right|\le C \left\| u-v\right\| _{L^{2} }^{a_{3} } \left\| u-v\right\| _{H^{s} }^{\gamma _{3} } \left(1+\left\| u\right\| _{H^{s} }^{\eta } +\left\| v\right\| _{H^{s} }^{\eta } \right),\] 
\noindent where the parameters satisfy the system
\[\left\{ {\begin{array}{*{20}{l}}
{{a_3} + {\gamma _3} \ge 2,\quad {a_3} \ge 0,0 < {\gamma _3} < 2,}\\
{{\gamma _3} + \eta  \le 2,\quad \,\,\eta  \ge 0.}
\end{array}} \right.\]

(\romannumeral4) For a sufficiently large positive constant $b$, the deterministic equation 
\begin{equation}\label{GrindEQ__3_3_}
D_t^\beta u =  - {\left( { - \Delta } \right)^s}u + b\Delta u + \zeta \left( u \right),\quad u\left( {x,0} \right) = {u_0}
\end{equation}
has a unique global solution $u$ with trajectories in ${L^2}\left( {0,T;{H^S}\left( {{\mathbb T^d}} \right)} \right) \cap C\left( {\left[ {0,T} \right];{L^2}\left( {{\mathbb T^d}} \right)} \right)$, and the ${L^2}$-norm of the solution $u$ is less than infinity when the initial data $u\left( {x,0} \right) = {u_0}$ is bounded, closed and convex. 
\end{hypothesis}

We state that $C$ in the inequality of conditions (\romannumeral1)-(\romannumeral3) denotes the constant independent of the other parameters, which is uniformly denoted by the symbol $C$ to avoid abuse of the symbols.

\begin{remark}\label{3.2}
By combining hypothesis (\romannumeral2), Lemma \ref{lemma2.6}  and Young inequality, we can conclude that exist some real constants $Q >0$, $a'=2a_{2} /\left(2-\gamma _{2} \right)>0$, any solutions of \eqref{GrindEQ__3_3_} hold
    \begin{gather*}
        \begin{array}{l} {D_{t}^{\beta } \left\| u\right\| _{L^{2} }^{2} =\int _{{\mathbb T}^{d} }D_{t}^{\beta } u^{2} dx \le \int _{{\mathbb T}^{d} }2u D_{t}^{\beta } u dx} \\{\rm  \; \; \; \; \; \; \; \; \; \; \; \; \; \; \; \; =-2\left\| \left(-\Delta \right)^{s/2} u\right\| _{L^{2} }^{2} -2b\left\| \nabla u\right\| _{L^{2} }^{2} +2\left\langle \zeta\left(u\right),u\right\rangle } \\ {{\rm   \; \; \; \; \; \; \; \; \; \; \; \; \; \; \; \; }\lesssim  -2\left\| \left(-\Delta \right)^{s/2} u\right\| _{L^{2} }^{2} -2b\left\| \nabla u\right\| _{L^{2} }^{2} +2\left(1+\left\| u\right\| _{H^{s} }^{\gamma _{2} } \right)\left(1+\left\| u\right\| _{L^{2} }^{a_{2} } \right)}. \end{array}
    \end{gather*}
Combining the inequality \eqref{GrindEQ__3_60_}(see \cite{flandoli2021delayed}),
\begin{equation}\label{GrindEQ__3_60_}
\left\| u\right\| _{H^{s} }^{2} \le 2^{s-1} \left(\left\| u\right\| _{L^{2} }^{2} +\left\| \left(-\Delta \right)^{s/2} u\right\| _{L^{2} }^{2} \right)
\end{equation}
and Poincar{\'e} inequality that we take the constants as $\frac{1}{Z}$, namely
\begin{equation}\label{GrindEQ__3_61_}
\left\| u-\int _{{\mathbb T}^{3} }u\left({x,t}\right)dx \right\| _{L^{2} }^{2} \le \frac{1}{Z} \left\| \nabla u\right\| _{L^{2} }^{2},
\end{equation}
the above inequality can be simplified as
\begin{equation}\label{GrindEQ__3_5_}
\begin{array}{l}
D_t^\beta \left\| u \right\|_{{L^2}}^2 \lesssim  - 2\left\| {{{\left( { - \Delta } \right)}^{s/2}}u} \right\|_{{L^2}}^2 - 2b\left\| {\nabla u} \right\|_{{L^2}}^2 + {2^{2 - s}}\left\| u \right\|_{{H^s}}^2 + Q\left( {1 + \left\| u \right\|_{{L^2}}^{a'}} \right)\\
{\rm{\; \; \; \; \; \; \; \; \; \;\; \; \; \; \;}} \le  - 2\left( {2{Z}b - 1} \right)\left\| u \right\|_{{L^2}}^2 + Q\left( {1 + \left\| u \right\|_{{L^2}}^{a'}} \right).
\end{array}
\end{equation}
\end{remark}

On the basis of the above, we discuss the model after adding multiplicative noise to equation \eqref{GrindEQ__3_50_}, which takes the form of
\begin{equation}\label{GrindEQ__3_51_}
D_t^\beta u = \left[ { - {{\left( { - \Delta } \right)}^s}u + \zeta \left( u \right)} \right]dt + A\sum\limits_{m \in Z_0^d} {\sum\limits_{j = 1}^{d - 1} {{\theta _m}{\sigma _{m,j}}\left( x \right)} }  \circ dW_t^{m,j},
\end{equation}
where $A$ is taken as $\sqrt {\frac{d}{{\left( {d - 1} \right)\left\| \theta  \right\|_{{\ell ^2}}^2}} \cdot b}$ in order to simplify the following proof, the exact meaning of its value is shown in the proof in the Section \ref{sec:4}. The parameters $\beta,\;s,\;\zeta\; $ are defined in the same way as model \eqref{GrindEQ__3_50_}. Next, Next, we focus on the remaining parameters in the equation \eqref{GrindEQ__3_51_} and the conditions that the parameters satisfy, which are similar to those in \cite{flandoli2021delayed,flandoli2021high}.

Denoting positive constants by $\mathbb Z_ + ^d $, correspondingly, we denote negative constants by $\mathbb Z_ - ^d $. For  $m \in \mathbb Z_ + ^d$, let $\theta  = \left\{ { \cdots ,{\theta _{ - m}}, \cdots ,{\theta _m}, \cdots } \right\}$ be a square addable sequence that satisfies ${\left\| \theta  \right\|_{{\ell ^2}}} < \infty$ and ${\theta _m} = {\theta _{ - m}}$. The vector fields ${\sigma _{m,j}}$ are defined as ${\sigma _{m,j}} = {q_{m,j}}{e^{2\pi mi \cdot x}}$. To sure that ${\sigma _{m,j}}$ are divergence free, we take \footnote{The choice of vector fields and complex Brownian motion is guided by Section 2.2 in \cite{galeati2020convergence}, and based on the method used in Section 2.3 for converting the Stratonovich form to the It{\^o} form, the equation is further simplified by choosing parameters $A$ of a particular form.}
\[\left\{ {\begin{array}{*{20}{l}}
{{q_{m,j}} \cdot m = 0,\quad \,m \in Z_ + ^d,}\\
{{q_{m,j}} = {q_{ - m,j}},\quad m \in Z_ - ^d,}
\end{array}} \right.\]
where $j = 1,2, \cdots ,d - 1.$ Summarizing the above, for any $x \in {{{\mathbb R^d}} \mathord{\left/
 {\vphantom {{{\mathbb R^d}} {{\mathbb Z^d}}}} \right.
 \kern-\nulldelimiterspace} {{\mathbb Z^d}}}$, we have 
\begin{equation}\label{GrindEQ__3_70_}
\begin{array}{l}
\sum\limits_m {\sum\limits_{j=1}^{d - 1} {\theta _m^2\left( {{\sigma _{m,j}}\left( x \right) \otimes {\sigma _{ - m,j}}\left( x \right)} \right)} } \, = \sum\limits_m {\theta _m^2\sum\limits_{j=1}^{d - 1} {\left( {{q_{m,j}}{e^{2\pi im \cdot x}} \otimes {q_{ - m,j}}{e^{ - 2\pi im \cdot x}}} \right)} } \\{\rm{ \;\;\;\;\;\;}}
 = \sum\limits_m {\theta _m^2\left( {{q_{m,1}} \otimes {q_{ - m,1}} +  \cdots  + {q_{m,d - 1}} \otimes {q_{ - m,d - 1}}} \right)}  = \frac{{d - 1}}{d}\left\| \theta  \right\|_{{\ell ^2}}^2{I_d}.{\rm{ }}
\end{array}
 \end{equation}
Assume $\left( {\Omega ,\cal F,\mathbb P} \right)$ is a probability space with filtration $\left\{ {{\mathcal{F}_t}:t \ge 0} \right\}$, ${W_t^{m,j}}$ is the complex Brownian motion defined on this probability space \cite{flandoli2021high}, namely
   \[{W_t^{m,j}} =\left\{\begin{array}{c} {B_t^{m,j} +iB_t^{-m,j} ,\; m\in {\rm {\mathbb Z}}_{+}^{d} ;} \\ {B_t^{-m,j} -iB_t^{m,j} ,\; m\in {\rm {\mathbb Z}}_{-}^{d} .} \end{array}\right. \]
Hence it clearly holds $\bar W_t^{m,j} = W_t^{ - m,j}$. Furthermore, we assume that ${W^{{m_1},{j_1}}},{W^{{m_2},{j_2}}}$ are independent whenever ${j_1} \ne {j_2}$ or ${m_1} \ne  - {m_2}$, namely, we have
\[{\left[ {{W^{{m_1},{j_1}}},{W^{{m_2},{j_2}}}} \right]_t} = 2t{\delta _{{j_1} - {j_2}}}{\delta _{{m_1} + {m_2}}}.\]

Based on the relationship between the Stratonovich and It{\^o} integrals \cite{galeati2020convergence}, we can convert the equation \eqref{GrindEQ__3_51_} to an equivalent It{\^o} form, namely
\[\begin{array}{l}
D_t^\beta u = \left[ { - {{\left( { - \Delta } \right)}^s}u + \zeta \left( u \right)} \right]dt + A\sum\limits_{m \in Z_0^d} {\sum\limits_{j = 1}^{d - 1} {{\theta _m}{\sigma _{m,j}}\left( x \right)} } \nabla u \circ dW_t^{m,j}\\
 = \left[ { - {{\left( { - \Delta } \right)}^s}u + \zeta \left( u \right)} \right]dt + A\sum\limits_{m \in Z_0^d} {\sum\limits_{j = 1}^{d - 1} {{\theta _m}{\sigma _{m,j}}\left( x \right)} } \nabla u \cdot dW_t^{m,j} + \frac{A}{2}\sum\limits_{m \in Z_0^d} {\sum\limits_{j = 1}^{d - 1} {{\theta _m}d\left[ {{\sigma _{m,j}}\left( x \right)\nabla u,W_t^{m,j}} \right]} } \\
 = \left[ { - {{\left( { - \Delta } \right)}^s}u + \zeta \left( u \right) + b\Delta u} \right]dt + A\sum\limits_{m \in Z_0^d} {\sum\limits_{j = 1}^{d - 1} {{\theta _m}{\sigma _{m,j}}\left( x \right)} } \nabla u \cdot dW_t^{m,j}
\end{array}\]
The simplification of the last step holds due to the condition \eqref{GrindEQ__3_70_} and the value of the constant $A$.

In the following, we concentrate on the equation
\begin{equation}\label{GrindEQ__3_71_}
D_t^\beta u = \left[ { - {{\left( { - \Delta } \right)}^s}u + \zeta \left( u \right) + b\Delta u} \right]dt + A\sum\limits_{m \in Z_0^d} {\sum\limits_{j = 1}^{d - 1} {{\theta _m}{\sigma _{m,j}}\left( x \right)} } \nabla u \cdot dW_t^{m,j}
\end{equation}
with initial value is ${u_0}$. To prevent confusion, denote the solution of equation \eqref{GrindEQ__3_71_} by ${u_t}\left( {{u_0},\theta } \right)$ and the solution of equation \eqref{GrindEQ__3_3_} by ${u_t}\left( {{u_0}} \right)$. Here are the two main conclusions of this paper.

\begin{thm}\label{theorem3.1} For bounded, closed and convex initial data $u\left( {x,0} \right) = {u_0}$, there exists a square summable real sequence ${\left\{ {{\theta _m}} \right\}_m}$, such that the equation 
\[D_t^\beta u = \left[ { - {{\left( { - \Delta } \right)}^s}u + \zeta \left( u \right)} \right]dt + A\sum\limits_{m \in Z_0^d} {\sum\limits_{j = 1}^{d - 1} {{\theta _m}{\sigma _{m,j}}\left( x \right)} }  \circ dW_t^{m,j}\]
have a unique solution, where $\zeta$ is an unknown nonlinear function satisfies Hypothesis \ref{hypothesis3.1}. Moreover, for any $0 < T < \infty $, we can find positive parameters $S$ and $b$ big enough, the solution satisfies
\begin{equation}
\mathbb P\left( {{{\left\| {u_t}\left( {{u_0},\theta } \right) \right\|}_{C\left( {\left[ {0,T} \right];{H^{ - \gamma }}} \right)}} < S} \right) \ge 1 - \varepsilon,\quad \forall t \in \left[ {0,T} \right],\quad \forall \varepsilon  \in \left( {0,1} \right). 
\end{equation}
\end{thm}

In conjunction with Theorem \ref{theorem3.1}, we give further assumptions to derive the other main conclusion of this paper, namely the delayed effect of the perturbation of the noise on the explosion time of the solution to the deterministic equation. Suppose the special sequence ${\left\{ {{\theta _m}} \right\}_m}$ taken in the above theorem is ${\left\{ {{\theta ^S}} \right\}_S}$.

\begin{thm}\label{theorem3.2} For bounded, closed and convex initial data $u\left( {x,0} \right) = {u_0}$, $\zeta$ satisfies Hypothesis \ref{hypothesis3.1}, the non-negative parameter $b$ can be taken large enough. Assume there exist constants $\lambda$ and $K$ greater than 0 and sufficiently large, such that the${L^2}$-norm of the solution of the time fractional deterministic equation
\[D_t^\beta u =  - {\left( { - \Delta } \right)^s}u + b\Delta u + \zeta \left( u \right)\]
holds the inequality
\[{\left\| {{u_t}\left( {{u_0}} \right)} \right\|_{{L^2}}} \le K{\left\| {{u_0}} \right\|_{{L^2}}}{e^{ - \lambda t}}.\]
And the solution of the time fractional stochastic differential equation
\begin{equation}\label{1}
D_t^\beta u = \left[ { - {{\left( { - \Delta } \right)}^s}u + \zeta \left( u \right)} \right]dt + A\sum\limits_{m \in Z_0^d} {\sum\limits_{j = 1}^{d - 1} {{\theta ^S}{\sigma _{m,j}}\left( x \right)} }  \circ dW_t^{m,j}
\end{equation}
exists a pathwise unique global solution for sufficiently small initial data ${u_0}$.  Then, with great probability, the lifespan of the pathwise unique solution of the Eqs. \eqref{1} can be taken to infinity.

\end{thm}

\section{Examples}
\label{sec:4}
Two examples of nonlinear fractional-order models are given in this section to justify Hypothesis \ref{hypothesis3.1}, and in this section, we take $d=3$.

\subsection{Fractional Keller-Segel Equation}

The Keller-Segel model is the partial differential equation proposed by Keller and Segel to model chemotaxis in cellular systems , which is also one of the most important partial differential systems for learning about chemotactic aggregation \cite{keller1970initiation}.
The classical expression of the model takes the form
\begin{equation}\label{4.1}
\left\{ {\begin{array}{*{20}{l}}
{{\partial _t}\rho  = \Delta \rho  - \chi \nabla  \cdot \left( {\rho \nabla c} \right){\rm{\; \;  }}in\;\Omega ,}\\
{{\partial _t}c = \gamma \Delta c + \beta \rho  - \mu c,}\\
{\rho \left( {0,\; \cdot \;} \right) = {\rho _0} \ge 0,\quad c\left( {0,\; \cdot \;} \right) = {c_0} \ge 0,}\\
{{\partial _x}\rho \left( {t,\; \cdot \;} \right) = {\partial _x}c\left( {t,\; \cdot \;} \right) = 0{\rm{\; \;  }}on\;\partial \Omega ,}
\end{array}} \right.
\end{equation}
where $\chi ,\gamma ,\mu ,\beta $ are positive constants. For dimension $d=1$, the solutions of the system \eqref{4.1} are regular, while $d$ takes other values, the solutions with finite-time blow-up in the case of sufficiently large initial data, cf. \cite{escudero2006fractional,jager1992explosions,winkler2013finite}. Moreover,  Winkler in \cite{winkler2013finite} gives the set of initial values that make the solutions of the equation \eqref{4.1} blow up.

As an extension of model \eqref{4.1}, the fractional Keller-Segel model is suitable for situations where cell motility cannot be described by random walk \cite{escudero2006fractional}. For example, Carvalho et al. \cite{RN84} introduced the mild solution of the equation
\[\left\{ {\begin{array}{*{20}{l}}
{D_t^\beta u - v\Delta u + \left( {u \cdot \nabla } \right)u + \nabla p = f,\quad {\mathop{\rm in}\nolimits} \;{R^n},\;\;t > 0,}\\
{u\left( {x,0} \right) = {u_0}\left( x \right) \ge 0,\quad \quad {\mathop{\rm in}\nolimits} \;{R^n},\;\;t > 0}\\
{\nabla  \cdot u = 0,\quad \quad {\mathop{\rm in}\nolimits} \;{R^n},}
\end{array}} \right.\]
where $\beta  \in \left( {0,1} \right) $. Many properties of mild solutions, including existence, uniqueness, decay and regularity, have also been demonstrated. Next, we propose a class of fractional Keller-Segel equations and verify that Hypothesis \ref{hypothesis3.1} holds.

\begin{lemma}\label{lemma4.1} Hypothesis 2.1 holds for the time fractional Keller-Segel model
\begin{equation}\label{GrindEQ__4_1_}
\left\{ {\begin{array}{*{20}{l}}
{D_t^\beta \rho  = \Delta \rho  - \nabla \cdot\left( {\rho \nabla c} \right),}\\
{ - \Delta c = \rho  - \int_\Omega  \rho  \left( x \right)dx,}\\
{\nabla \cdot\rho  = 0,}
\end{array}} \right.
\end{equation}
where $\beta  \in \left( {0,1} \right)$,  the regular bounded domain $\Omega  = {{\mathbb T}^3}$.
\end{lemma}

\noindent $\mathbf{Proof.}$  Apparently, the first equation in the above system can be simplified as
\[\begin{array}{l}
D_t^\beta \rho  = \Delta \rho  - \nabla  \cdot \left( {\rho \nabla c} \right) = \Delta \rho  - \nabla  \cdot \left[ {\rho \nabla  \cdot \left( { - {\Delta ^{ - 1}}} \right)\left( {\rho  - {\rho _\Omega }} \right)} \right]\\{\rm \; \; \; \; \; \; \;}
 = \Delta \rho  - \nabla  \cdot \left[ {\rho {\nabla ^{ - 1}}\left( {\rho  - {\rho _\Omega }} \right)} \right].
\end{array}\]
which is equivalent to taking $s = 1$, $\zeta \left( \rho  \right) =  - \nabla  \cdot \left[ {\rho {\nabla ^{ - 1}}\left( {\rho  - {\rho _\Omega }} \right)} \right]$ in equation \eqref{GrindEQ__3_50_}, note that ${\rho _\Omega }$ denotes $\int_{{\mathbb T^3}} {\rho \left( {x,t} \right)dx}$. To begin with, we show that the nonlinear mapping $\zeta $ in this model satisfies the conditions (\romannumeral1)-(\romannumeral3). From the Sobolev embedding inequality, we obtain the inclusion relation between different Sobolev spaces as follows.
\begin{equation}\label{2}
{H^1}\left( {{{\mathbb T}^3}} \right) \hookrightarrow {H^{3/4}}\left( {{{\mathbb T}^3}} \right) \hookrightarrow {L^4}\left( {{{\mathbb T}^3}} \right),\quad H^{7/4} \left({\rm {\mathbb T}}^{3} \right)\hookrightarrow L^{\infty } \left({\rm {\mathbb T}}^{3} \right),\quad H^{1/2} \left({\rm {\mathbb T}}^{3} \right)\hookrightarrow L^{3} \left({\rm {\mathbb T}}^{3} \right).
\end{equation}
By considering the properties of the norm, we deduce that
\[\begin{array}{l}
{\left\| {\zeta \left( {{\rho _1}} \right) - \zeta \left( {{\rho _2}} \right)} \right\|_{{H^{ - 1}}}} = {\left\| { - \nabla  \cdot \left[ {{\rho _1}{\nabla ^{ - 1}}\left( {{\rho _1} - {\rho _{2\Omega }}} \right)} \right] + \nabla  \cdot \left[ {{\rho _2}{\nabla ^{ - 1}}\left( {{\rho _2} - {\rho _{2\Omega }}} \right)} \right]} \right\|_{{H^{ - 1}}}}\\{\rm \;\;\;\;\;\;\;\;\;\;\;\;\;\;\;\;\;\;\;\;\;\;}
 = {\left\| {\nabla  \cdot \left( {{\rho _1}{\nabla ^{ - 1}}{\rho _1}} \right) - \nabla  \cdot \left[ {{\rho _1}{\nabla ^{ - 1}}{\rho _{1\Omega }}} \right] - \nabla  \cdot \left( {{\rho _2}{\nabla ^{ - 1}}{\rho _2}} \right) + \nabla  \cdot \left[ {{\rho _2}{\nabla ^{ - 1}}{\rho _{2\Omega }}} \right]} \right\|_{{H^{ - 1}}}}\\{\rm \;\;\;\;\;\;\;\;\;\;\;\;\;\;\;\;\;\;\;\;\;\;}
 \le {\left\| {\nabla  \cdot \left( {{\rho _1}{\nabla ^{ - 1}}{\rho _1}} \right) - \nabla  \cdot \left( {{\rho _2}{\nabla ^{ - 1}}{\rho _2}} \right)} \right\|_{{H^{ - 1}}}} + {\left\| {\nabla  \cdot \left[ {{\rho _1}{\nabla ^{ - 1}}{\rho _{1\Omega }}} \right] - \nabla  \cdot \left[ {{\rho _2}{\nabla ^{ - 1}}{\rho _{2\Omega }}} \right]} \right\|_{{H^{ - 1}}}}\\{\rm \;\;\;\;\;\;\;\;\;\;\;\;\;\;\;\;\;\;\;\;\;\;}
 \le {\left\| {\left( {{\rho _1}{\nabla ^{ - 1}}{\rho _1}} \right) - \left( {{\rho _2}{\nabla ^{ - 1}}{\rho _2}} \right)} \right\|_{{L^2}}} + {\left\| {{\rho _1}{\nabla ^{ - 1}}{\rho _{1\Omega }} - {\rho _2}{\nabla ^{ - 1}}{\rho _{2\Omega }}} \right\|_{{L^2}}}.
\end{array}\]
where ${\rho _{1\Omega }}\left( t \right) = {\rho _{1\Omega }}\left( 0 \right) = \int_\Omega  {{\rho _1}\left(t, x \right)dx}$ and $ {\rho _{2\Omega }}\left( t \right) = {\rho _{2\Omega }}\left( 0 \right) = \int_\Omega  {{\rho _2}\left(t, x \right)dx}$ are both constants. Next we scale the norm with the triangle inequality and results \eqref{2}, the following formula holds \footnote{Decomposition and scaling of the norm by means of a method similar to that of Lemma 2.1 in \cite{flandoli2021delayed}. }
\[\begin{array}{l}
{\left\| {\left( {{\rho _1}{\nabla ^{ - 1}}{\rho _1}} \right) - \left( {{\rho _2}{\nabla ^{ - 1}}{\rho _2}} \right)} \right\|_{{L^2}}} = {\left\| {{\rho _1}{\nabla ^{ - 1}}{\rho _1} - {\rho _2}{\nabla ^{ - 1}}{\rho _1} + {\rho _2}{\nabla ^{ - 1}}{\rho _1} - {\rho _2}{\nabla ^{ - 1}}{\rho _2}} \right\|_{{L^2}}}\\{\rm \;\;\;\;\;\;\;\;\;\;\;\;\;\;\;\;\;\;\;\;\;\;\;\;\;\;\;\;\;\;\;\;\;\;\;\;\;\;\;\;\;\;\;\;\;}
 \lesssim {\left\| {{\rho _1} - {\rho _2}} \right\|_{{L^2}}}{\left\| {{\nabla ^{ - 1}}{\rho _1}} \right\|_{{L^\infty }}} + {\left\| {{\rho _2}} \right\|_{{L^4}}}{\left\| {{\nabla ^{ - 1}}\left( {{\rho _1} - {\rho _2}} \right)} \right\|_{{L^4}}}\\{\rm \;\;\;\;\;\;\;\;\;\;\;\;\;\;\;\;\;\;\;\;\;\;\;\;\;\;\;\;\;\;\;\;\;\;\;\;\;\;\;\;\;\;\;\;\;}
 \lesssim {\left\| {{\rho _1} - {\rho _2}} \right\|_{{L^2}}}{\left\| {{\rho _1}} \right\|_{{H^{{3 \mathord{\left/
 {\vphantom {3 4}} \right.
 \kern-\nulldelimiterspace} 4}}}}} + {\left\| {{\rho _2}} \right\|_{{H^{{3 \mathord{\left/
 {\vphantom {3 4}} \right.
 \kern-\nulldelimiterspace} 4}}}}}{\left\| {\left( {{\rho _1} - {\rho _2}} \right)} \right\|_{{L^2}}}.
\end{array}\]
Synthesizing the above narrative, the inequality
\[{\left\| {\zeta \left( {{\rho _1}} \right) - \zeta\left( {{\rho _2}} \right)} \right\|_{{H^{ - 1}}}} \lesssim {\left\| {{\rho _1} - {\rho _2}} \right\|_{{L^2}}}\left( {1 + {{\left\| {{\rho _1}} \right\|}_{{H^{3/4}}}} + {{\left\| {{\rho _2}} \right\|}_{{H^{3/4}}}}} \right).\]
holds.  To verify the condition (\romannumeral1), we let ${\rho _2} = 0$ in the above equation, namely
\[{\left\| {\zeta \left( {{\rho _1}} \right)} \right\|_{{H^{ - 1}}}} \lesssim {\left\| {{\rho _1}} \right\|_{{L^2}}}\left( {1 + {{\left\| {{\rho _1}} \right\|}_{{H^{3/4}}}}} \right) \lesssim \left( {1 + {{\left\| {{\rho _1}} \right\|}_{{L^2}}}} \right)\left( {1 + {{\left\| {{\rho _1}} \right\|}_{{H^{1}}}}} \right),\]
which means that (\romannumeral1) holds with ${a_1} = 1,{\gamma _1} = 1/4.$

For the proof of condition (\romannumeral2), we need to use the key inequality  $\left\langle {f,{\nabla ^{ - 1}}f\cdot\nabla f} \right\rangle  =  - \frac{1}{2}\left\langle {{f^2},\nabla \cdot\left( {{\nabla ^{ - 1}}f} \right)} \right\rangle $ and the nature of inner product, it holds
\[\begin{array}{l}
\left| {\left\langle {\zeta \left( \rho  \right),\rho } \right\rangle } \right| = \left| {\left\langle { - \nabla  \cdot \left[ {\rho {\nabla ^{ - 1}}\left( {\rho  - {\rho _\Omega }} \right)} \right],\rho } \right\rangle } \right| \\{\rm \; \; \; \; \;\; \; \; \; \;\; \; \; \; \;}
 = \left| {\left\langle { - \nabla \rho  \cdot {\nabla ^{ - 1}}\left( {\rho  - {\rho _\Omega }} \right),\rho } \right\rangle  + \left\langle { - \rho \left( {\rho  - {\rho _\Omega }} \right),\rho } \right\rangle } \right|\\{\rm \; \; \; \; \;\; \; \; \; \;\; \; \; \; \;}
 \lesssim \left| {\int_\Omega  {{\rho ^2}\left( x \right)\left( {\rho \left( x \right) - {\rho _\Omega }\left( x \right)} \right)dx} } \right|\\{\rm \; \; \; \; \;\; \; \; \; \;\; \; \; \; \;}
\lesssim \left\| \rho  \right\|_{{L^3}}^3.
\end{array}\]
We can obtain from the embedding theorem that the norm satisfies the formula
\begin{equation}\label{6}
\left\| \rho  \right\|_{{H^{1/2}}}^3 \lesssim\left\| \rho  \right\|_{{L^2}}^{3/2}\left\| \rho  \right\|_{{H^{1/2}}}^{3/2}  \lesssim \left\| \rho  \right\|_{{L^2}}^{3/2}\left\| \rho  \right\|_{{H^1}}^{3/2},
\end{equation}
therefore it clearly holds
\[\left| {\left\langle {\zeta \left( \rho  \right),\rho } \right\rangle } \right| \lesssim \left( {1 + \left\| \rho  \right\|_{{L^2}}^{3/2}} \right)\left( {1 + \left\| \rho  \right\|_{{H^1}}^{3/2}} \right),\]
which means that (\romannumeral2) holds with ${a_2} = {\gamma _2} = 3/2.$

It is clear that
\[{\left| {\left\langle {\zeta \left( {{\rho _1}} \right) - \zeta \left( {{\rho _2}} \right),{\rho _1} - {\rho _2}} \right\rangle } \right| \le {{\left\| {\zeta \left( {{\rho _1}} \right) - \zeta \left( {{\rho _2}} \right)} \right\|}_{{H^{ - 1}}}}{{\left\| {{\rho _1} - {\rho _2}} \right\|}_{{H^1}}}},\]
then we can use the result of condition (\romannumeral1) to get the inequality easily
\[{\left| {\left\langle {\zeta \left( {{\rho _1}} \right) - \zeta \left( {{\rho _2}} \right),{\rho _1} - {\rho _2}} \right\rangle } \right| \lesssim {{\left\| {{\rho _1} - {\rho _2}} \right\|}_{{L^2}}}{{\left\| {{\rho _1} - {\rho _2}} \right\|}_{{H^1}}}\left( {1 + {{\left\| {{\rho _1}} \right\|}_{{H^1}}} + {{\left\| {{\rho _2}} \right\|}_{{H^1}}}} \right)},\]
which means that (\romannumeral3) is satisfied with $a_{3} =\gamma _{3} =\eta =1.$ 

Substituting $\rho \left( t \right) - {\rho _{{{\mathbb T}^3}}}\left( t \right) = \rho \left( t \right) - m$ into the equation, we can obtain the same equation of the form \ref{GrindEQ__3_3_}, then from the Remark \ref{3.2}, it follows that
\[D_t^\beta \left\| {\rho  - m} \right\|_{{L^2}}^2 \le  - 2\left( {4Zb - 1} \right)\left\| {\rho  - m} \right\|_{{L^2}}^2 + Q\left[ {1 + {{\left( {\left\| {\rho  - m} \right\|_{{L^2}}^2} \right)}^{{{a'} \mathord{\left/
 {\vphantom {{a'} 2}} \right.
 \kern-\nulldelimiterspace} 2}}}} \right].\]
We set $D_t^\beta f\left( t \right) =  - 2\left( {4Zb - 1} \right)f\left( t \right) + Q\left[ {1 + x{{\left( t \right)}^{{{a'} \mathord{\left/
 {\vphantom {{a'} 2}} \right.
 \kern-\nulldelimiterspace} 2}}}} \right]$, by using Property \ref{property2.1} then the following identity holds
\[f\left( t \right) = f\left( 0 \right) + \frac{1}{{\Gamma \left( \beta  \right)}}\int_0^t {{{\left( {t - s} \right)}^{\beta  - 1}}\left[ {Q\left( {1 + f{{\left( s \right)}^{{{a'} \mathord{\left/
 {\vphantom {{a'} 2}} \right.
 \kern-\nulldelimiterspace} 2}}}} \right) - 2\left( {4Zb - 1} \right)x\left( s \right)} \right]ds} .\]
Assume that $f\left( 0 \right) = \left\| {\rho \left( 0 \right) - {\rho _{{{\mathbb T}^3}}}\left( 0 \right)} \right\|_{{L^2}}^2$ is finite, ${{\rho _{{{\mathbb T}^3}}}\left( 0 \right)} \in \mathbb R$, we can take $b$ sufficiently large such that $f\left( t \right) < \infty$. As an application of Lemma \ref{lemma2.4}, (\romannumeral4) holds.

\subsection{Fractional Fisher-KPP Equation}

The fisher-kpp model is a typical example of reaction-diffusion equation, which which is named after Fisher \cite{fisher1937wave}, Kolmogorov, and Petrovsky and Piskunov \cite{tikhomirov1991study}. The Fisher-KPP model has an extremely wide range of applications, such as the study on the spatial and temporal spread of epidemics, the spatial spread of invasive species, cf.\cite{RN92,henry2000fractional}. From \cite{angenent1995asymptotic}, we know that a typical form of such equation is
\[\frac{\partial u\left(x,t\right)}{\partial t} =D\frac{\partial ^{2} u\left(x,t\right)}{\partial x^{2} } +ru\left(x,t\right)\left(1-u\left(x,t\right)\right),\; \quad D>0,r>0.\]

With the rapid development of related theories, classical Fisher-KPP models have many limitations in solving practical problems, so scholars have introduced fractional order into such equations, which are considered to be more effective for settling sub-diffusion problems \cite{metzler2000random}. There are also a lot of research results on the explosion of solutions of such equations. For example, 
Ahmad et al. \cite{ahmad2015time} respectively discussed under what initial value conditions does the solution blow up in finite time, and under what initial value conditions does the solution satisfy bounded and global existence.  Xu et al. \cite{xu2018extremely} showed the blow-up phenomenon and the conditions for its appearance, moreover, they stated that the lower the order, the earlier the blow-up comes. Next, we consider the fractional Fisher-KPP equation in 3D case.

\begin{lemma}\label{lemma4.2} Hypothesis 3.1 holds for the time fractional order model
\begin{equation} \label{GrindEQ__44_} 
\left\{ {\begin{array}{*{20}{l}}
{D_t^\beta u = \Delta u + {u^2} - u,}\\
{\int_{{\mathbb T^3}} {u\left( {x,0} \right)} dx \in \left( { - \infty ,1} \right),\quad {{\left\| {u\left( {x,0} \right) - \int_{{\mathbb T^3}} {u\left( {x,0} \right)dx} } \right\|}_{{L^2}}} \le \sqrt {{\delta _0}} ,}
\end{array}} \right.
\end{equation}
where $0 < \beta  < 1$, $0 \le {\delta _0} < \infty$.
\end{lemma}

\noindent $\mathbf{Proof.}$ 
The system \eqref{GrindEQ__44_} is equivalent to taking $s = 1$, $\zeta \left( \rho  \right) = {u^2} - u$ in equation \eqref{GrindEQ__3_50_}. Let us first list a few important inequalities that are used in the proof below. The interpolation inequality (see \cite{shou1986interpolation})
\begin{equation}\label{GrindEQ__4_10_}
{\left\| u \right\|_{{L^2}}} \le \left\| u \right\|_{{H^\alpha }}^{\frac{\beta }{{\alpha  + \beta }}}\left\| u \right\|_{{H^{ - \beta }}}^{\frac{\alpha }{{\alpha  + \beta }}},
\end{equation} 
\noindent the Sobolev inequality (see \cite{guo2022existence})
\begin{equation}\label{GrindEQ__4_9_}
\left\| u \right\|_{{L^{{p_k}}}\left( {{\mathbb R^n}} \right)}^{{p_k}} \le \left\| u \right\|_{{L^{\frac{{n{p_k}}}{{n - 2}}}}\left( {{\mathbb R^n}} \right)}^{{\varsigma _1}{p_k}}\left\| u \right\|_{{L^{{p_k} - 1}}\left( {{\mathbb R^n}} \right)}^{{p_k}\left( {1 - {\varsigma _1}} \right)},
\end{equation}
with ${\varsigma _1} = {{\left( {\frac{1}{{{p_{k - 1}}}} - \frac{1}{{{p_k}}}} \right)} \mathord{\left/{\vphantom {{\left( {\frac{1}{{{p_{k - 1}}}} - \frac{1}{{{p_k}}}} \right)} {\left( {\frac{1}{{{p_{k - 1}}}} - \frac{{n - 2}}{{n{p_k}}}} \right)}}} \right.
 \kern-\nulldelimiterspace} {\left( {\frac{1}{{{p_{k - 1}}}} - \frac{{n - 2}}{{n{p_k}}}} \right)}} \sim O\left( 1 \right),\quad 1 - {\varsigma _1} \sim O\left( 1 \right),\quad {p_k} \to \infty$,\quad and the H{\"o}lder's inequality (see \cite{RN87})
 \[\int_a^b {\left| {f\left( x \right)g\left( x \right)} \right|dx}  \le {\left( {\int_a^b {{{\left| {f\left( x \right)} \right|}^p}dx} } \right)^{\frac{1}{p}}}{\left( {\int_a^b {{{\left| {g\left( x \right)} \right|}^q}dx} } \right)^{\frac{1}{q}}},\]
where $p,q > 0$ and $\frac{1}{p} + \frac{1}{q} = 1$.

The proof of the (\romannumeral1)-(\romannumeral3) of Hypothesis \ref{hypothesis3.1} are similar to that of \cite{flandoli2021delayed}. Using the embedding relationship  ${H^{1/2}}\left( {{{\mathbb T}^3}} \right) \hookrightarrow {L^3}\left( {{{\mathbb T}^3}} \right)$, we deduce 
  \[\begin{array}{l}
{\left\| {\zeta \left( {{u_1}} \right) - \zeta \left( {{u_2}} \right)} \right\|_{{H^{ - 1}}}} = {\left\| {\left( {u_1^2 - {u_1}} \right) - \left( {u_2^2 - {u_2}} \right)} \right\|_{{H^{ - 1}}}}\\{\rm \;\;\;\;\;\;\;\;\;\;\;\;\; \;\;\;\;\;\;\;\;\;\;\;\;\;\;\;\;\;\;\;\;}
 \lesssim {\left\| {\left( {u_1^2 - u_2^2} \right) - \left( {{u_1} - {u_2}} \right)} \right\|_{{H^{  -{1 \mathord{\left/
 {\vphantom {1 2}} \right.
 \kern-\nulldelimiterspace} 2}}}}} \\{\rm \;\;\;\;\;\;\;\;\;\;\;\;\; \;\;\;\;\;\;\;\;\;\;\;\;\;\;\;\;\;\;\;\;}
 \lesssim {\left\| {{u_1} + {u_2}} \right\|_{{H^{  {1 \mathord{\left/
 {\vphantom {1 2}} \right.
 \kern-\nulldelimiterspace} 2}}}}}{\left\| {{u_1} - {u_2}} \right\|_{{H^{  {1 \mathord{\left/
 {\vphantom {1 2}} \right.
 \kern-\nulldelimiterspace} 2}}}}} + {\left\| {{u_1} - {u_2}} \right\|_{{H^{  {1 \mathord{\left/
 {\vphantom {1 2}} \right.
 \kern-\nulldelimiterspace} 2}}}}}\\{\rm \;\;\;\;\;\;\;\;\;\;\;\;\; \;\;\;\;\;\;\;\;\;\;\;\;\;\;\;\;\;\;\;\;}
 \lesssim {\left\| {{u_1} - {u_2}} \right\|_{{H^{{1 \mathord{\left/
 {\vphantom {1 2}} \right.
 \kern-\nulldelimiterspace} 2}}}}}\left( {1 + {{\left\| {{u_1}} \right\|}_{{H^{{1 \mathord{\left/
 {\vphantom {1 2}} \right.
 \kern-\nulldelimiterspace} 2}}}}} + {{\left\| {{u_2}} \right\|}_{{H^{{1 \mathord{\left/
 {\vphantom {1 2}} \right.
 \kern-\nulldelimiterspace} 2}}}}}} \right).
\end{array}\]
Then we take ${u_2} = 0$, namely
\[{\left\| {\zeta \left( {{u_1}} \right)} \right\|_{{H^{ - 1}}}} \lesssim {\left\| {{u_1}} \right\|_{{H^{{1 \mathord{\left/
 {\vphantom {1 2}} \right.
 \kern-\nulldelimiterspace} 2}}}}}\left( {1 + {{\left\| {{u_1}} \right\|}_{{H^{{1 \mathord{\left/
 {\vphantom {1 2}} \right.
 \kern-\nulldelimiterspace} 2}}}}}} \right),\]
this inequality not only shows that the mapping from $H^{1/2}$ to $ H^{-1}$ has continuity, but also indicates that (\romannumeral1) holds with $a_{1} =1,\gamma _{1} =1/2.$

With respect to hypothesis (\romannumeral2), by using the inequality \eqref{6}, we can easily conclude 
\[\begin{array}{l}
\left| {\left\langle {\zeta \left( u \right),u} \right\rangle } \right| = \left| {\left\langle {{u^2} - u,u} \right\rangle } \right| = \int_{{\mathbb T^3}} {\left( {{u^2} - u} \right)udx} \\{\rm \;\;}\quad \quad \quad \quad  \le \left\| u \right\|_{{L^3}}^3 \lesssim \left\| u \right\|_{{H^1}}^{3/2}\left\| u \right\|_{{L^2}}^{3/2}
\end{array}\]
which implies condition (\romannumeral2) is satisfied with
${a_2} = 3/2,{\gamma _2} = 3/2$.

Then we test assumption (\romannumeral3), since
\[\begin{array}{l}
\left| {\left\langle {\zeta \left( {{u_1}} \right) - \zeta \left( {{u_2}} \right),{u_1} - {u_2}} \right\rangle } \right| = \left| {\left\langle {\left( {u_1^2 - {u_1}} \right) - \left( {u_2^2 - {u_2}} \right),{u_1} - {u_2}} \right\rangle } \right| \\{\rm \;\;\;\;\;\;\;\;\;\;\;\;\; \;\;\;\;\;\;\;\;\;\;\;\;\;\;\;\;\;\;\;\;\;\;\;\;\;\;\;\;}
= \left| {\left\langle {\left( {u_1^2 - u_2^2} \right) - \left( {{u_1} - {u_2}} \right),{u_1} - {u_2}} \right\rangle } \right|\\{\rm \;\;\;\;\;\;\;\;\;\;\;\;\; \;\;\;\;\;\;\;\;\;\;\;\;\;\;\;\;\;\;\;\;\;\;\;\;\;\;\;\;}
\lesssim {\left\| {{u_1} + {u_2}} \right\|_{{L^3}}}\left\| {{u_1} - {u_2}} \right\|_{{L^3}}^2\\{\rm \;\;\;\;\;\;\;\;\;\;\;\;\; \;\;\;\;\;\;\;\;\;\;\;\;\;\;\;\;\;\;\;\;\;\;\;\;\;\;\;\;}
 \lesssim {\left\| {{u_1} + {u_2}} \right\|_{{H^{1/2}}}}\left\| {{u_1} - {u_2}} \right\|_{{H^{{1 \mathord{\left/
 {\vphantom {1 2}} \right.
 \kern-\nulldelimiterspace} 2}}}}^2\\{\rm \;\;\;\;\;\;\;\;\;\;\;\;\; \;\;\;\;\;\;\;\;\;\;\;\;\;\;\;\;\;\;\;\;\;\;\;\;\;\;\;\;}
 \lesssim \left( {{{\left\| {{u_1}} \right\|}_{{H^1}}} + {{\left\| {{u_2}} \right\|}_{{H^1}}}} \right){\left\| {{u_1} - {u_2}} \right\|_{{H^1}}}{\left\| {{u_1} - {u_2}} \right\|_{{L^2}}},
\end{array}\]
which implies (\romannumeral3) is satisfied with $a_{3} =1,\gamma _{3} =1, \eta =1.$

Finally, we verify that (\romannumeral4) holds. For equation  $D_t^\beta u = \Delta u + b\Delta u + {u^2} - u$  over ${\mathbb T}^{3}$, we integrate both sides over ${\mathbb T}^{3}$ simultaneously, 
\[D_t^\beta {u_{{{\mathbb T}^3}}} = \left( {1 + b} \right)\Delta {u_{{{\mathbb T}^3}}} + \left\| u \right\|_{{L^2}}^2 - {u_{{{\mathbb T}^3}}} = \left\| u \right\|_{{L^2}}^2 - u_{{{\mathbb T}^3}}^2 + u_{{{\mathbb T}^3}}^2 - {u_{{{\mathbb T}^3}}} = \left\| {u - {u_{{{\mathbb T}^3}}}} \right\|_{{L^2}}^2 + u_{{{\mathbb T}^3}}^2 - {u_{{{\mathbb T}^3}}}.\]
Subtracting the two equations, we say
\[\begin{array}{l}
D_t^\beta \left( {u - {u_{{{\mathbb T}^3}}}} \right) = \left( {1 + b} \right)\Delta \left( {u - {u_{{{\mathbb T}^3}}}} \right) + {u^2} - u - \left( {\left\| u \right\|_{{L^2}}^2 - {u_{{{\mathbb T}^3}}}} \right)\\{\rm \;\;\;\;\;\;\;\;}
 = \left( {1 + b} \right)\Delta \left( {u - {u_{{{\mathbb T}^3}}}} \right) + \left( {{u^2} - \left\| u \right\|_{{L^2}}^2} \right) - \left( {u - {u_{{{\mathbb T}^3}}}} \right),
\end{array}\]
since $\left\| {u - {u_{{\mathbb T^3}}}} \right\|_{{L^2}}^2 \ge 0$ and $\left\| u \right\|_{{L^2}}^2 \ge 0$ hold constantly, as an application of Lemma \ref{lemma2.6},  we can easily deduce
\[\begin{array}{l}
D_t^\beta \left\| {u - {u_{{{\mathbb T}^3}}}} \right\|_{{L^2}}^2 \le 2\int_{{\mathbb T^3}} {\left( {u - {u_{{\mathbb T^3}}}} \right)D_t^\beta \left( {u - {u_{{\mathbb T^3}}}} \right)dx} \\
{\rm{ \;\;\;\;\;\;\;\;\;\;\;\;\;\;\;\;\;\;\;\;\;               }}=- 2\left( {1 + b} \right)\left| {\left| {\nabla \left( {u - {u_{{{\mathbb T}^3}}}} \right)} \right|} \right|_{{L^2}}^2 - 2\left\| {u - {u_{{{\mathbb T}^3}}}} \right\|_{{L^2}}^2 + 2\int_{{{\mathbb T}^3}} {\left( {{u^2} - \left\| u \right\|_{{L^2}}^2} \right)\left( {u - {u_{{{\mathbb T}^3}}}} \right)dx} \\
{\rm{ \;\;\;\;\;\;\;\;\;\;\;\;\;\;\;\;\;\;\;\;\;               }} \le  - 2\left( {1 + b} \right)\left| {\left| {\nabla \left( {u - {u_{{{\mathbb T}^3}}}} \right)} \right|} \right|_{{L^2}}^2 + 2\int_{{{\mathbb T}^3}} {{u^2}\left( {u - {u_{{{\mathbb T}^3}}}} \right)dx}.
\end{array}\]
Next, we work on the term $\int_{{{\mathbb T}^3}} {{u^2}\left( {u - {u_{{{\mathbb T}^3}}}} \right)dx}$, note that
\[\begin{array}{l}
{u^2}\left( {u - {u_{{{\mathbb T}^3}}}} \right) = {\left( {u - {u_{{{\mathbb T}^3}}}} \right)^3} + 2{u_{{{\mathbb T}^3}}}{u^2} - 3u_{{{\mathbb T}^3}}^2\left( {u - {u_{{{\mathbb T}^3}}}} \right) + u_{{{\mathbb T}^3}}^3\\{\rm \;\;\;\;\;\;\;\;\;\;\;\;\;\;\;\;\;\;}
={\left( {u - {u_{{{\mathbb T}^3}}}} \right)^3} + 2{u_{{{\mathbb T}^3}}}{\left( {u - {u_{{{\mathbb T}^3}}}} \right)^2} + u_{{{\mathbb T}^3}}^2\left( {u - {u_{{{\mathbb T}^3}}}} \right)
\end{array}\]
since $\int_{{{\mathbb T}^3}} {\left( {u - {u_{{{\mathbb T}^3}}}} \right)dx}  = 0$, then, combining \eqref{GrindEQ__4_10_} and \eqref{GrindEQ__4_9_}  yields \footnote{The proof that this inequality holds is inspired by Proposition 2.4 in \cite{flandoli2021delayed}, from which the following system can be obtained.}
\[\begin{array}{l}
\int_{{{\mathbb T}^3}} {{u^2}\left( {u - {u_{{{\mathbb T}^3}}}} \right)dx}  = \int_{{{\mathbb T}^3}} {{{\left( {u - {u_{{T^3}}}} \right)}^3}dx}  + 2{u_{{{\mathbb T}^3}}}\int_{{{\mathbb T}^3}} {{{\left( {u - {u_{{{\mathbb T}^3}}}} \right)}^2}dx} \\
{\rm{\;\;\;\;\;\;\;\;}} \le \left\| {u - {u_{{{\mathbb T}^3}}}} \right\|_{{L^3}}^3 + 2{u_{{{\mathbb T}^3}}}\left\| {u - {u_{{{\mathbb T}^3}}}} \right\|_{{L^2}}^2\\
{\rm{\;\;\;\;\;\;\;\;}} \le \left\| {u - {u_{{{\mathbb T}^3}}}} \right\|_{{H^{{1 \mathord{\left/
 {\vphantom {1 2}} \right.
 \kern-\nulldelimiterspace} 2}}}}^3 + 2{u_{{{\mathbb T}^3}}}\left\| {u - {u_{{T^3}}}} \right\|_{{L^2}}^2\\
{\rm{\;\;\;\;\;\;\;\;}} \lesssim \left\| {u - {u_{{{\mathbb T}^3}}}} \right\|_{{L^2}}^{3/2}\left\| {u - {u_{{{\mathbb T}^3}}}} \right\|_{{H^1}}^{3/2} + 2{u_{{{\mathbb T}^3}}}\left\| {u - {u_{{{\mathbb T}^3}}}} \right\|_{{L^2}}^2\\
{\rm{\;\;\;\;\;\;\;\;}} \le P\left\| {u - {u_{{{\mathbb T}^3}}}} \right\|_{{L^2}}^6 + \left\| {\nabla \left( {u - {u_{{{\mathbb T}^3}}}} \right)} \right\|_{{L^2}}^2 + 2{u_{{{\mathbb T}^3}}}\left\| {u - {u_{{{\mathbb T}^3}}}} \right\|_{{L^2}}^2,
\end{array}\]
where $P$ is an unimportant constant. Then, the above inequality can be transformed into
\[\begin{array}{l}
D_t^\beta \left\| {u - {u_{{{\mathbb T}^3}}}} \right\|_{{L^2}}^2 \le  - 2\left( {1 + b} \right)\left| {\left| {\nabla \left( {u - {u_{{{\mathbb T}^3}}}} \right)} \right|} \right|_{{L^2}}^2 + 2P\left\| {u - {u_{{{\mathbb T}^3}}}} \right\|_{{L^2}}^6 \\{\rm{\;\;\;\;\;\;\;\;\;\;\;\;\;\;\;\;\;\;\;\;\;\;\;\;\;\;}}
+ 2\left\| {\nabla \left( {u - {u_{{{\mathbb T}^3}}}} \right)} \right\|_{{L^2}}^2 + 4{u_{{{\mathbb T}^3}}}\left\| {u - {u_{{{\mathbb T}^3}}}} \right\|_{{L^2}}^2\\{\rm{\;\;\;\;\;\;\;\;\;\;\;\;\;\;\;\;\;\;\;\;\;\;\;\;}}
 =  - 2b\left| {\left| {\nabla \left( {u - {u_{{{\mathbb T}^3}}}} \right)} \right|} \right|_{{L^2}}^2 + 2P\left\| {u - {u_{{{\mathbb T}^3}}}} \right\|_{{L^2}}^6 + 4{u_{{{\mathbb T}^3}}}\left\| {u - {u_{{{\mathbb T}^3}}}} \right\|_{{L^2}}^2\\{\rm{\;\;\;\;\;\;\;\;\;\;\;\;\;\;\;\;\;\;\;\;\;\;\;\;}}
\le  - 2b\left| {\left| {\nabla u} \right|} \right|_{{L^2}}^2 + 2\left\| {u - {u_{{T^3}}}} \right\|_{{L^2}}^6 + 4{u_{{T^3}}}\left\| {u - {u_{{{\mathbb T}^3}}}} \right\|_{{L^2}}^2\\{\rm{\;\;\;\;\;\;\;\;\;\;\;\;\;\;\;\;\;\;\;\;\;\;\;\;}}
 \le \left( { - 2bZ + 4{u_{{{\mathbb T}^3}}}} \right)\left\| {u - {u_{{{\mathbb T}^3}}}} \right\|_{{L^2}}^2 + 2P\left\| {u - {u_{{{\mathbb T}^3}}}} \right\|_{{L^2}}^6,
\end{array}\]
the last step is obtained from \eqref{GrindEQ__3_61_}.
In summary, we can obtain the system of equations
\[\left\{ {\begin{array}{*{20}{l}}
{D_t^\beta {u_{{{\mathbb T}^3}}} = \left\| {u - {u_{{{\mathbb T}^3}}}} \right\|_{{L^2}}^2 + u_{{{\mathbb T}^3}}^2 - {u_{{{\mathbb T}^3}}},}\\
{D_t^\beta \left\| {u - {u_{{{\mathbb T}^3}}}} \right\|_{{L^2}}^2 \le \left( { - 2bZ + 4{u_{{{\mathbb T}^3}}}} \right)\left\| {u - {u_{{{\mathbb T}^3}}}} \right\|_{{L^2}}^2 + 2P\left\| {u - {u_{{{\mathbb T}^3}}}} \right\|_{{L^2}}^6.}
\end{array}} \right.\]
Since the proof methods of the two cases $0 \le \int_{{{\mathbb T}^3}} {{u_0}\left( x \right)dx}< 1$ and $\int_{{{\mathbb T}^3}} {{u_0}\left( x \right)dx}< 0$ are not the same, we shall discuss them separately.

For $0 \le \int_{{{\mathbb T}^3}} {{u_0}\left( x \right)dx}  \le {k_0} < 1$, we set 
\[{M_1} = \left\{ {\left. {t \ge 0} \right|0 \le {u_{{{\mathbb T}^3}}}\left( t \right) < {k_0} + \xi {\rm{ }}} \right\} \cap \left\{ {\left. {t \ge 0} \right|\left\| {u - {u_{{{\mathbb T}^3}}}} \right\|_{{L^2}}^2 < {\delta _0} + \xi } \right\},{\rm \;\;} \forall \xi  \in \left( {0,1} \right), \]
for any $t \in {M_1}$, it holds
\[D_t^\beta \left\| {u - {u_{{{\mathbb T}^3}}}} \right\|_{{L^2}}^2 \le \left[ { - 2bZ + 4\left( {{k_0} + \xi } \right) + 2P{{\left( {{\delta _0} + \xi } \right)}^2}} \right]\left\| {u - {u_{{{\mathbb T}^3}}}} \right\|_{{L^2}}^2.\]
Take constant $b$ large enough such that $2bZ - 4\left( {{k_0} + \xi } \right) - 2P\left( {{\delta _0} + \xi } \right) > 0$, we then use the Lemma \ref{lemma2.5}, yield
\[\left\| {u\left( {t,x} \right) - {u_{{{\mathbb T}^3}}}\left( t \right)} \right\|_{{L^2}}^2 \le \left\| {{u_0} - {u_{{{\mathbb T}^3}}}\left( 0 \right)} \right\|_{{L^2}}^2 \le {\delta _0},\]
hence, we say
\[D_t^\beta {u_{{{\mathbb T}^3}}} = \left\| {u - {u_{{{\mathbb T}^3}}}} \right\|_{{L^2}}^2 + u_{{{\mathbb T}^3}}^2 - {u_{{{\mathbb T}^3}}} \Rightarrow D_t^\beta {u_{{{\mathbb T}^3}}} + {u_{{{\mathbb T}^3}}} = \left\| {u - {u_{{{\mathbb T}^3}}}} \right\|_{{L^2}}^2 + u_{{{\mathbb T}^3}}^2 \le {\delta _0} + 1 + {e^{ - \lambda t}},\]
where $\lambda $ is a non-negative constant. Applying Lemma \ref{lemma2.5} and Property \ref{property2.1}, we get
\[\begin{array}{l}
{u_{{{\mathbb T}^3}}}\left( t \right) \le {u_{{{\mathbb T}^3}}}\left( 0 \right) + \frac{1}{{\Gamma \left( \beta  \right)}}\int_0^t {{{\left( {t - s} \right)}^{\beta  - 1}}\left( {{\delta _0} + 1 + {e^{ - \lambda s}}} \right)ds} \\
{\rm{\;\;\;\;\;\;\;\;\;\;\;      }} \le {u_{{{\mathbb T}^3}}}\left( 0 \right) + \frac{K}{{\Gamma \left( \beta  \right)}}\int_0^t {{{\left( {t - s} \right)}^{\beta  - 1}}{e^{ - \lambda s}}ds} \\
{\rm{\;\;\;\;\;\;\;\;\;\;\;     }} \le {k_0} + \xi  < 1,
\end{array}\]
the above inequality is valid because the growth rate of exponential function is faster than that of power function. It is possible to take $K$ and $\lambda$ that satisfy inequalities ${e^{ - \lambda s}} + 1 + {\delta _0} \le K{e^{ - \lambda s}}$ and $\frac{K}{{\Gamma \left( \beta  \right)}}\int_0^t {{{\left( {t - s} \right)}^{\beta  - 1}}{e^{ - \lambda s}}ds}  \le \xi$, which implies the conclusion holds.

Next we discuss the second case that ${u_{{{\mathbb T}^3}}} < 0$. Set 
\[{M_2} = \left\{ {\left. {t \ge 0} \right|{u_{{{\mathbb T}^3}}}\left( t \right) < 0{\rm{ }}} \right\} \cap \left\{ {\left. {t \ge 0} \right|\left\| {u - {u_{{{\mathbb T}^3}}}} \right\|_{{L^2}}^2 < {\delta _0} + \xi } \right\}, {\rm \;\;} \forall \xi  \in \left( {0,1} \right),\]
for any $t \in {M_2}$, the following inequality is valid,
\[D_t^\beta \left\| {u - {u_{{{\mathbb T}^3}}}} \right\|_{{L^2}}^2 \le \left[ { - 2bZ + 2P{{\left( {{\delta _0} + \xi } \right)}^2}} \right]\left\| {u - {u_{{{\mathbb T}^3}}}} \right\|_{{L^2}}^2.\]
As above, when $b$ is large enough such that the inequality  $2bZ - 2P{\left( {{\delta _0} + \xi } \right)^2} > 0$ holds, we can obtain
\[\left\| {u\left( {t,x} \right) - {u_{{{\mathbb T}^3}}}\left( t \right)} \right\|_{{L^2}}^2 \le \left\| {{u_0} - {u_{{{\mathbb T}^3}}}\left( 0 \right)} \right\|_{{L^2}}^2 \le {\delta _0}.\]

Combing the conclusion of the two cases above, it is clear that
\[\left\{ {\begin{array}{*{20}{l}}
{{u_{{{\mathbb T}^3}}}\left( t \right) \le 1,}\\
{\left\| {u\left( {t,x} \right) - {u_{{{\mathbb T}^3}}}\left( t \right)} \right\|_{{L^2}}^2 \le \left\| {{u_0} - {u_{{{\mathbb T}^3}}}\left( 0 \right)} \right\|_{{L^2}}^2 \le {\delta _0},}
\end{array}} \right.\]
therefore, we have 
\[\mathop {\sup }\limits_{0 \le t \le T} {\left\| {u\left( t \right)} \right\|_{{L^2}}} \le \mathop {\sup }\limits_{0 \le t \le T} {\left\| {u\left( {t,x} \right) - {u_{{T^3}}}\left( t \right)} \right\|_{{L^2}}} + \mathop {\sup }\limits_{0 \le t \le T} {\left\| {{u_{{T^3}}}\left( t \right)} \right\|_{{L^2}}}< \infty.\]
the deduction of this inequality stems from the application of the triangle inequality, and also shows that condition (\romannumeral4) holds.\qed

\begin{remark}
The proof of Lemma \ref{lemma4.2} is based on the verification in \cite{flandoli2021delayed}, the main feature of this paper is applying the triangle inequality, Cordoba-Cordoba inequality and Lemma \ref{lemma2.5} to scale the ${L^2}$-norm of the solution for the equation \eqref{GrindEQ__44_} to obtain the conclusion.  
\end{remark}

In the above discussion, we restricted $\int_{{{\mathbb T}^3}} {{u_0}\left( x \right)dx}$ to the range less than 1. The following proposition will explain why the value $\int_{{{\mathbb T}^3}} {{u_0}\left( x \right)dx} \ge 1$ cannot be taken.

\begin{proposition}\label{proposition4.2} The time fractional equation
\begin{equation}\label{GrindEQ__4_13_}
\left\{ {\begin{array}{*{20}{l}}
{D_t^\beta u\left( t \right) = \left( {1 + b} \right)\Delta u\left( t \right) + {u^2}\left( t \right) - u\left( t \right),}\\
{{u_{{\mathbb T^3}}}\left( 0 \right) = \int_{{\mathbb T^3}} {u\left( {x,0} \right)dx}  > 1,}
\end{array}} \right.
\end{equation}
experience an explosion in a finite time.
\end{proposition}

\noindent $\mathbf{Proof.}$ The first equation in System \eqref{GrindEQ__4_13_} is integrated over ${\rm {\mathbb T}}^{3} $ on both the left-hand and right-hand sides at the same time, we have
\[D_{t}^{\beta } u_{{\mathbb T}^{3} } =(1+b)\Delta u_{{\mathbb T}^{3} } +\int _{_{{\mathbb T}^{3} } }u^{2} dx -u_{{\mathbb T}^{3} },\] 
recall the Jensen inequality (see \cite{dragomir1994some}), we have the conclusion
\[\int_{{\mathbb T^3}} {{u^2}\left( {t,x} \right)dx}  \ge {\left( {\int_{{\mathbb T^3}} {u\left( {t,x} \right)dx} } \right)^2} = u_{{\mathbb T^3}}^2\left( t \right).\]
Note that $\Delta u_{{\mathbb T}^{3} } =\int _{{\mathbb T}^{3} }\Delta u\left(x,t\right)dx =0$,  combining the above conclusions, we know from the Theorem 3.2 in \cite{RN76} that for initial values satisfying the condition $x\left( 0 \right) > 1$, the solution of the equation 
\begin{equation}\label{42}
D_t^\beta x\left( t \right) = {x^2}\left( t \right) - x\left( t \right)
\end{equation}
satisfies the identity
\[\mathop {\lim }\limits_{t \to {T_1}} x\left( t \right) =  + \infty,\]
where ${T_1}$ is a constant less than infinity, which implies that the blow-up of the solution of Eq. \eqref{42} occurs in a limited time. This also means that the blow-up of the solution of Eq. \eqref{42} occurs in a limited time. 

We set ${x\left( 0 \right) = \int_{{\mathbb T^3}} {u\left( {x,0} \right)dx}}$,  as an application of Lemma \ref{lemma2.4}, it holds
\[\mathop {\lim }\limits_{t \to {T_1}} {\left\| {u\left( t \right)} \right\|_{{L^2}}} \ge \mathop {\lim }\limits_{t \to {T_1}} {u_{{\mathbb T^3}}}\left( t \right) \ge \mathop {\lim }\limits_{t \to {T_1}} x\left( t \right) =  + \infty ,\]
namely, the conclusion holds. \qed

\section{Existence and uniqueness of solution}
\label{sec:5}

The emphasis of this study is on the effect of the existence of noise on the blow-up time to the deterministic equations. In this section, we shall use Galerkin approximation and a priori estimates methods to prove the existence and uniqueness of solution to the fractional stochastic differential equations.

While the function $\zeta$ in the equation is nonlinear, its solution may not exist globally but locally \cite{flandoli2021delayed}. In order to discuss the global solution of the equation, we introduce a smooth non-increasing cut-off function $L_{S} $, which has the form

\begin{equation}\label{GrindEQ__5_1}
{L_S}\left( {{{\left\| u \right\|}_{{H^{ - \gamma }}}}} \right) = \left\{ {\begin{array}{*{20}{l}}
{1,\quad \quad 0 \le {{\left\| u \right\|}_{{H^{ - \gamma }}}} \le S,}\\
{0,\quad \quad {{\left\| u \right\|}_{{H^{ - \gamma }}}} > S + 1,}
\end{array}} \right.
\end{equation}
note that $L_{S} \left(\left\| u\right\| _{H^{-\gamma } } \right)$ is equivalent to $L_{S} \left(u\right)$ when take $\gamma $  to be sufficiently small.
Next, we work on fractional order stochastic differential equations with cut-off function $L_{S} $, the system can be rewritten as 
\begin{equation} \label{GrindEQ__5_2_} 
\left\{\begin{array}{l} {D_{t}^{\beta } u={\left[ { - {{\left( { - \Delta } \right)}^s}u + b\Delta u + {L_S}\left( {{{\left\| u \right\|}_{{H^{ - \gamma }}}}} \right)\zeta \left( u \right)} \right]dt + A\sum\limits_m {\sum\limits_{j=1}^{d - 1} {{\theta _m}{\sigma _{m,j}}\nabla u\left( t \right)dW_t^{m,j}} } } ,} \\ u\left( {x,0} \right) = u\left( 0 \right). \end{array}\right.  
\end{equation} 
where $s \ge 1,\frac{1}{2} <\beta <1$, $A = \sqrt {\frac{d}{{\left( {d - 1} \right)\left\| \theta  \right\|_{{\ell ^2}}^2}} \cdot b} $, $b>0$. In the next expression, we write $\Lambda =\left(-\Delta \right)^{1/2} $.

\begin{definition}\label{definition5.1}
A filtered probability space $\left( {\Omega ,{\rm{{\cal F}}},\left( {{{\rm{{\cal F}}}_t}} \right),\mathbb P} \right)$ is considered, $\left\{ {{W^{m,j}}:m \in Z_0^d,j = 1, \cdots ,d - 1} \right\}$ stand for  a collection of independent complex Brownian motions, which is defined on $\Omega $. The $\left( {{{\mathcal{F}}_t}} \right)$-progressively measurable process $u$ is said as a solution to the equation \eqref{GrindEQ__5_2_}, with trajectories of class ${L^\infty }\left( {0,T;{L^2}\left( {{{\mathbb T}^d}} \right)} \right) \cap {L^2}\left( {0,T;{H^s}\left( {{{\mathbb T}^d}} \right)} \right)$, $\mathbb P - a.s$ , 
\[\begin{array}{*{20}{l}}
{\left\langle {u\left( t \right),\varphi } \right\rangle  = \left\langle {{u_0},\varphi } \right\rangle  - \frac{1}{{\Gamma \left( \beta  \right)}}\int_0^t {{{\left( {t - \tau } \right)}^{\beta  - 1}}\left[ {\left\langle {{\Lambda ^s}u\left( \tau  \right),{\Lambda ^s}\varphi } \right\rangle  + b\left\langle {\nabla u\left( \tau  \right),\nabla \varphi } \right\rangle } \right]d\tau } }\\{\rm \;\;\;\;\;\;\;\;\;\;\;\;\;\;}
{ + \frac{1}{{\Gamma \left( \beta  \right)}}\int_0^t {{{\left( {t - \tau } \right)}^{\beta  - 1}}{L_S}\left( {{{\left\| u \right\|}_{{H^{ - \gamma }}}}} \right)\left\langle {\zeta \left( {u\left( \tau  \right)} \right),\varphi } \right\rangle } d\tau }\\{\rm \;\;\;\;\;\;\;\;\;\;\;\;\;\;}
{ - \frac{A}{{\Gamma \left( \beta  \right)}}\sum\limits_m {\sum\limits_{j = 1}^{d - 1} {{\theta _m}\int_0^t {{{\left( {t - \tau } \right)}^{\beta  - 1}}\left\langle {u\left( \tau  \right),{\sigma _{m,j}}\cdot\nabla \varphi } \right\rangle } dW_\tau ^{m,j}} } }
\end{array}\]
for all $t \in \left[ {0,T} \right]$, where ${u_0} \in {L^2}\left( {{\mathbb T^d}} \right)$ and the function $\varphi  \in {H^s}\left( {{\mathbb T^d}} \right)$ is divergence free.
\end{definition}

\begin{remark}\label{remark5.1} 
To facilitate the following discussion, we briefly note
\[G\left( t \right) = \frac{1}{{\Gamma \left( \beta  \right)}}\int_0^t {{{\left( {t - \tau } \right)}^{\beta  - 1}}\left[ { - {\Lambda ^{2s}}u\left( \tau  \right) + b\Delta u\left( \tau  \right) + {L_S}\left( {{{\left\| {u\left( \tau  \right)} \right\|}_{{H^{ - \gamma }}}}} \right)\zeta \left( {u\left( \tau  \right)} \right)} \right]} d\tau, \]
\[M\left( t \right) = \frac{A}{{\Gamma \left( \beta  \right)}}\int_0^t {\sum\limits_m {\sum\limits_{j = 1}^{d-1} {{\theta _m}{\sigma _{m,j}}{{\left( {t - \tau } \right)}^{\beta  - 1}} \cdot \nabla u\left( \tau  \right)dW_\tau ^{m,j}} } }. \]
\end{remark}
The problem is first solved on a finite dimensional space, namely employing the Galerkin approximations. Set ${{\rm H}_M} = {\mathop{\rm span}\nolimits} \left\{ {{e^{2\pi mi \cdot x}}:\left| m \right| \le M} \right\}$, which is a finite dimensional space and is also the subspace of ${L^2}\left( {{{\mathbb T}^d}} \right)$. Then, define the orthogonal projection ${\Pi _M}:{\rm H} \to {{\rm H}_M}$, it can be  be expressed as \[{\Pi _M}h = \sum\limits_{m \in Z_0^d} {{h_m}{e^{2\pi mi \cdot x}}} ,\quad \quad \forall h \in {\rm H}.\]
From the above definition, it is clear that the for any $h \in {{\rm H}_M}$, ${\left\| h \right\|_{{L^2}}} = {\left\| h \right\|_{{H^s}}}$.  Then, we discuss the fractional order stochastic differential equations on ${{\rm H}_M}$.

\begin{lemma}\label{lemma5.1} The time fractional order stochastic differential equation  
\begin{equation}\label{GrindEQ__5_71_}
\left\{ {\begin{array}{*{20}{l}}
\begin{array}{l}
D_t^\beta {u_M}\left( t \right) = \left[ { - {\Lambda ^{2s}}{u_M}\left( t \right) + b\Delta {u_M}\left( t \right) + {L_S}\left( {{{\left\| {{u_M}\left( t \right)} \right\|}_{{H^{ - \gamma }}}}} \right){\Pi _M}\zeta \left( {{u_M}\left( t \right)} \right)} \right]dt\\
\quad \quad \quad \quad  - A\sum\limits_m {\sum\limits_{j = 1}^{d - 1} {{\theta _m}{\Pi _M}\left( {{\sigma _{m,j}}\cdot\nabla {u_M}\left( t \right)} \right)dW_t^{m,j},} } 
\end{array}\\
{{u_M}\left( 0 \right) = {\Pi _M}{u_0},}
\end{array}} \right.
\end{equation}
admits a unique solution, note that  $\theta \in \ell ^{2}$ and $u_{0} \in L^{2} \left({\rm {\mathbb T}}^{d} \right)$. Moreover, it holds the inequality
\begin{equation} \label{GrindEQ__5_11_} 
\mathop {\sup }\limits_{t \in \left[ {0,T} \right]} \left\| {u_M\left( t \right)} \right\|_{{L^2}}^2 \le K\left( {\left\| {{u_0}} \right\|_{{L^2}}^2,T,\gamma } \right)
\end{equation} 
with probability no less than $1 - \varepsilon $, the constant $ K\left( {\left\| {{u_0}} \right\|_{{L^2}}^2,T,\gamma } \right)$ indicates that it is related to the parameters ${\left\| {{u_0}} \right\|_{{L^2}}^2,T,\gamma }$ and is independent of the time $t$.
\end{lemma}

\noindent $\mathbf{Proof.}$ In the following expression, we shall denote the solution of Eq.\eqref{GrindEQ__5_71_} by ${{u_M}}$. By using Lemma \ref{lemma2.5}, we have
\[\begin{array}{l}
D_t^\beta \left\| {{u_M}} \right\|_{{L^2}}^2 \le 2\int_{{\mathbb T^d}} {{u_M}\left( t \right)D_t^\beta {u_M}\left( t \right)dx} =  - 2\left\| {{\Lambda ^s}{u_M}\left( t \right)} \right\|_{{L^2}}^2 - 2b\left\| {\nabla {u_M}\left( t \right)} \right\|_{{L^2}}^2 \\{\rm{   \;\;\;\;\;\;\;\;\;\;\;\;\;\;\;\;\;\;}}
+ 2{L_S}\left( {{{\left\| {{u_M}\left( t \right)} \right\|}_{{H^{ - \gamma }}}}} \right)\left\langle {\zeta \left( {{u_M}\left( t \right)} \right),{u_M}\left( t \right)} \right\rangle \\
{\rm{   \;\;\;\;\;\;\;\;\;\;\;\;\;\;\;\;\;\;}} - 2A\sum\limits_m {\sum\limits_{j=1}^{d - 1} {{\theta _m}} \left\langle {{u_M}\left( t \right),{\Pi _M}\left( {{\sigma _{m,j}} \cdot \nabla {u_M}\left( t \right)} \right)} \right\rangle dW_t^{m,j}}.
\end{array}\]
\noindent Since the projection $\Pi _{M}$ is orthogonal, and ${\sigma _{m,j}}$ is divergence free by the definition in Section \ref{sec:3}, we state 
\[\left\langle u_{M} \left(t\right),\Pi _{M} \left(\sigma _{m,j} \cdot \nabla u_{M} \left(t\right)\right)\right\rangle =\left\langle u_{M} \left(t\right),\sigma _{m,j} \cdot \nabla u_{M} \left(t\right)\right\rangle =0,\]  
considering the condition $b>0$, the above inequality can be simplified as
\begin{equation} \label{GrindEQ__5_13_} 
D_{t}^{\beta } \left\| u_{M} \right\| _{L^{2} }^{2} \le -2\left\| \Lambda ^{s} u_{M} \left(t\right)\right\| _{L^{2} }^{2} +2L_{S} \left(u_{M} \left(t\right)\right)\left\langle \zeta\left(u_{M} \left(t\right)\right),u_{M} \left(t\right)\right\rangle  
\end{equation} 

Next, we mainly work on the second term of \eqref{GrindEQ__5_13_}. Applying the inequality in the Hypothesis \ref{hypothesis3.1}  (\romannumeral2), we state 
\[\begin{array}{l}
2{L_S}\left( {{{\left\| {{u_M}\left( t \right)} \right\|}_{{H^{ - \gamma }}}}} \right)\left\langle {\zeta \left( {{u_M}\left( t \right)} \right),{u_M}\left( t \right)} \right\rangle  \lesssim 2{L_S}\left( {{{\left\| {{u_M}\left( t \right)} \right\|}_{{H^{ - \gamma }}}}} \right)\left( {1 + \left\| {{u_M}\left( t \right)} \right\|_{{H^s}}^{{\gamma _2} + {a_2}\frac{\gamma }{{s + \gamma }}}} \right)\\
\quad \quad \quad \quad \quad \quad \quad \quad \quad \quad \quad \quad \quad \quad \quad \quad \quad \times \left( {1 + \left\| {{u_M}\left( t \right)} \right\|_{{H^{ - \gamma }}}^{{a_2}\frac{s}{{s + \gamma }}}} \right)\\
\quad \quad \quad \quad \quad \quad \quad \quad \quad \quad \quad \quad \quad \quad \quad \quad \quad  \lesssim 2\left( {1 + {S^{{a_2}\frac{s}{{s + \gamma }}}}} \right)\left( {1 + \left\| {{u_M}\left( t \right)} \right\|_{{H^s}}^{{\gamma _2} + {a_2}\frac{\gamma }{{s + \gamma }}}} \right),
\end{array}\]
combining Young's inequality and \eqref{GrindEQ__3_60_}, the above equation can be further simplified as 
\[\begin{array}{l}
2{L_S}\left( {{{\left\| {{u_M}\left( t \right)} \right\|}_{{H^{ - \gamma }}}}} \right)\left\langle {\zeta \left( {{u_M}\left( t \right)} \right),{u_M}\left( t \right)} \right\rangle \; \le 2\left[ {1 + {S^{{a_2}\frac{s}{{s + \gamma }}}} + C{{\left( {1 + {S^{{a_2}\frac{s}{{s + \gamma }}}}} \right)}^{\frac{2}{{2 - {\gamma _2} - {a_2}\frac{\gamma }{{s + \gamma }}}}}}} \right]\\{\rm{ \;\;\;\;\;\;\;\;\;\;\;\;\;\;\;\;\;\;\;\;\;\;\;\;\;\;\;\;\;\;\;\;\;\;\;\;\;\;\;\;\;\;\;\;\;\;\;\;\;\;\;\;\;\;\;\;\;\;\;\;\;\;\;\;}}
 + \left\| {{u_M}\left( t \right)} \right\|_{{L^2}}^2 + \left\| {{\Lambda ^s}{u_M}\left( t \right)} \right\|_{{L^2}}^2.
\end{array}\]

For simplifying the proof, we abbreviate $2\left[ {1 + {S^{{a_2}\frac{s}{{s + \gamma }}}} + C{{\left( {1 + {S^{{a_2}\frac{s}{{s + \gamma }}}}} \right)}^{\frac{2}{{2 - {\gamma _2} - {a_2}\frac{\gamma }{{s + \gamma }}}}}}} \right]$  as the constant $\tilde C$, it is easy to obtain 
\[D_t^\beta \left\| {{u_M}} \right\|_{{L^2}}^2  + \left\| {{\Lambda ^s}{u_M}\left( t \right)} \right\|_{{L^2}}^2 \le  \left\| {{u_M}\left( t \right)} \right\|_{{L^2}}^2 + \tilde C.\]
Since it is obvious that $\left\| {{\Lambda ^s}{u_M}\left( \tau  \right)} \right\|_{{L^2}}^2 $ is greater than 0, we say
\[D_t^\beta \left\| {{u_M}} \right\|_{{L^2}}^2 \le  \left\| {{u_M}\left( t \right)} \right\|_{{L^2}}^2 + \tilde C.\]
Recall  Property \ref{property2.1}, yields
\[\left\| {{u_M}\left( t \right)} \right\|_{{L^2}}^2 \le \frac{1}{{\Gamma \left( \beta  \right)}}\int_0^t {{{\left( {t - \tau } \right)}^{\beta  - 1}}\left\| {{u_M}\left( \tau  \right)} \right\|_{{L^2}}^2d\tau }  + C'{t^\beta } + \left\| {u\left( 0 \right)} \right\|_{{L^2}}^2\]
where $\beta  \in \left( {1/2,1} \right)$ and ${C'>0}$. We apply the Theorem 2.4 in \cite{zhu2018new}, the inequality \eqref{GrindEQ__5_11_} holds.

We prove the pathwise uniqueness of the solution using  priori estimate. Suppose that $u_{1}$ and $u_{2}$ are two solutions of Eq.\eqref{GrindEQ__5_11_}, which correspond to the same Brownian motion and initial data ${u_0} \in {L^2}\left( {{\mathbb T^d}} \right)$. Set $\xi=u_{1} -u_{2}$, the following identity holds,
\[\begin{array}{*{20}{l}}
{\left\langle {\xi\left( t \right),\varphi } \right\rangle  =  - \frac{1}{{\Gamma \left( \beta  \right)}}\int_0^t {{{\left( {t - \tau } \right)}^{\beta  - 1}}} \left\langle {{\Lambda ^s}\xi \left( \tau  \right),{\Lambda ^s}\varphi } \right\rangle d\tau  + \frac{b}{{\Gamma \left( \beta  \right)}}\int_0^t {{{\left( {t - \tau } \right)}^{\beta  - 1}}} \left\langle {\nabla \xi u\left( \tau  \right),\nabla \varphi } \right\rangle d\tau }\\
{\;\;\;\;\;\;\;\;\;\;\;\;\; + \frac{1}{{\Gamma \left( \beta  \right)}}\int_0^t {{{\left( {t - \tau } \right)}^{\beta  - 1}}} \left\langle {{L_S}\left( {{{\left\| {{u_1}} \right\|}_{{H^{ - \gamma }}}}} \right)\zeta \left( {{u_1}\left( \tau  \right)} \right) - {L_S}\left( {{{\left\| {{u_2}} \right\|}_{{H^{ - \gamma }}}}} \right)\zeta \left( {{u_2}\left( t \right)} \right),\varphi } \right\rangle d\tau }\\
{\;\;\;\;\;\;\;\;\;\;\;\;\; - \frac{A}{{\Gamma \left( \beta  \right)}}\sum\limits_m {\sum\limits_{j = 1}^{d - 1} {{\theta _m}\int_0^t {{{\left( {t - \tau } \right)}^{\beta  - 1}}\left\langle {\xi\left( \tau  \right),{\sigma _{m,j}}\cdot\nabla \varphi } \right\rangle } dW_\tau ^{m,j}} } ,}
\end{array}\]
for any  $\varphi \in H^{s} \left({\rm {\mathbb T}}^{2} \right)$. Similar to \eqref{GrindEQ__5_13_}, it holds
\begin{equation}\label{GrindEQ__5_17_}
D_t^\beta \left\| {\xi \left( t \right)} \right\|_{{L^2}}^2 \le  - 2\left\| {{\Lambda ^s}\xi \left( t \right)} \right\|_{{L^2}}^2 + 2\left\langle {{L_S}\left( {{{\left\| {{u_1}} \right\|}_{{H^{ - \gamma }}}}} \right)\zeta \left( {{u_1}\left( t \right)} \right) - {L_S}\left( {{{\left\| {{u_2}} \right\|}_{{H^{ - \gamma }}}}} \right)\zeta \left( {{u_2}\left( t \right)} \right),\xi \left( t \right)} \right\rangle .
\end{equation}
Note that the second term in the Eq. \eqref{GrindEQ__5_17_} can be split into
\[\begin{array}{l}
{L_S}\left( {{{\left\| {{u_1}} \right\|}_{{H^{ - \gamma }}}}} \right)\zeta \left( {{u_1}\left( t \right)} \right) - {L_S}\left( {{{\left\| {{u_2}} \right\|}_{{H^{ - \gamma }}}}} \right)\zeta \left( {{u_2}\left( t \right)} \right) = {L_S}\left( {{{\left\| {{u_1}} \right\|}_{{H^{ - \gamma }}}}} \right)\zeta \left( {{u_1}\left( t \right)} \right) - {L_S}\left( {{{\left\| {{u_2}} \right\|}_{{H^{ - \gamma }}}}} \right)\zeta \left( {{u_1}\left( t \right)} \right)\\{\rm{   \;\;\;\;\;\;\;\;\;\;\;\;\;\;\;\;\;\;\;\;\;\;\;\;\;\;\;}}
 + {L_S}\left( {{{\left\| {{u_2}} \right\|}_{{H^{ - \gamma }}}}} \right)\zeta \left( {{u_1}\left( t \right)} \right) - {L_S}\left( {{{\left\| {{u_2}} \right\|}_{{H^{ - \gamma }}}}} \right)\zeta \left( {{u_2}\left( t \right)} \right)\\{\rm{   \;\;\;\;\;\;\;\;\;\;\;\;\;\;\;\;\;\;\;\;\;\;\;\;\;\;\;}}
 = \left[ {{L_S}\left( {{{\left\| {{u_1}} \right\|}_{{H^{ - \gamma }}}}} \right) - {L_S}\left( {{{\left\| {{u_2}} \right\|}_{{H^{ - \gamma }}}}} \right)} \right]\zeta \left( {{u_1}\left( t \right)} \right) + {L_S}\left( {{{\left\| {{u_2}} \right\|}_{{H^{ - \gamma }}}}} \right)\left[ {\zeta \left( {{u_1}\left( t \right)} \right) - \zeta \left( {{u_2}\left( t \right)} \right)} \right],
\end{array}\]
we shall discuss them separately.

Then we use Lagrange mean value theorem, hypothesis (\romannumeral1) and the condition that the ${L^2}$-norm of ${u_i}\left( {i = 1,2} \right)$ is bounded, we deduce
\[\begin{array}{*{20}{l}}
\begin{array}{l}
\left| {{L_S}\left( {{{\left\| {{u_1}} \right\|}_{{H^{ - \gamma }}}}} \right) - {L_S}\left( {{{\left\| {{u_2}} \right\|}_{{H^{ - \gamma }}}}} \right)} \right|\left| {\left\langle {\zeta \left( {{u_1}\left( t \right)} \right),\xi \left( t \right)} \right\rangle } \right| \mathbin\lesssim {\left\| {{{L'}_S}} \right\|_\infty }\left| {{{\left\| {{u_1}} \right\|}_{{H^{ - \gamma }}}} - {{\left\| {{u_2}} \right\|}_{{H^{ - \gamma }}}}} \right|\\
\;\;\;\;\;\;\;\;\;\;\;\;\;\;\;\;\;\;\;\;\;\;\;\;\;\;\;\;\;\;\;\;\;\;\;\;\;\;\;\;\;\;\;\;\;\;\;\;\;\;\;\;\;\;\;\;\;\;\;\;\;\;\;\;\;\;\;\;\;\;\;\;\;\;\;\;\;\;\;\; \times {\left\| {\zeta \left( {{u_1}\left( t \right)} \right)} \right\|_{{H^{ - s}}}}{\left\| {\xi \left( t \right)} \right\|_{{H^s}}}
\end{array}\\
{\;\;\;\;\;\;\;\;\;\;\;\;\;\;\;\;\;\;\;\;\;\;\;\;\;\;\;\;\;\;\;\;\;\;\;\;\;\;\;\;\;\;\;\;\;\;\;\;\;\;\;\;\;\;\;\;\;\;\;\;\;\;\;\;\;\;\;\;\;\;\;\;\;\;\;\;\;\; \mathbin\lesssim {{\left\| \xi  \right\|}_{{L^2}}}\left( {1 + {{\left\| {{u_1}\left( t \right)} \right\|}_{{H^s}}}} \right){{\left\| {\xi \left( t \right)} \right\|}_{{H^s}}}.}
\end{array}\]
Combining the inequality \eqref{GrindEQ__3_60_}, the above inequality can be deformed as
\begin{equation}\label{GrindEQ__5_18_}
\left| {{L_S}\left( {{{\left\| {{u_1}} \right\|}_{{H^{ - \gamma }}}}} \right) - {L_S}\left( {{{\left\| {{u_2}} \right\|}_{{H^{ - \gamma }}}}} \right)} \right|{\left\langle {\zeta \left( {{u_1}\left( t \right)} \right),\xi \left( t \right)} \right\rangle }{\le K\left( {1 + \left\| {{u_1}\left( t \right)} \right\|_{{H^s}}^2} \right)\left\| {\xi \left( t \right)} \right\|_{{L^2}}^2 + \frac{1}{2}\left\| {{\Lambda ^s}\xi \left( t \right)} \right\|_{{L^2}}^2.}
\end{equation}
For the another part, we apply the same inequalities as above and hypothesis (\romannumeral3), thus having 
\begin{equation} \label{GrindEQ__5_19_} 
\begin{array}{*{20}{l}}
{{L_S}\left( {{u_2}\left( t \right)} \right)\left| {\left\langle {\zeta \left( {{u_1}\left( t \right)} \right) - \zeta \left( {{u_2}\left( t \right)} \right),\xi \left( t \right)} \right\rangle } \right|}{\lesssim \left\| {\xi \left( t \right)} \right\|_{{H^s}}^{{\gamma _3}}\left\| {\xi\left( t \right)} \right\|_{{L^2}}^{{a_3}}\left( {1 + \left\| {{u_1}\left( t \right)} \right\|_{{H^s}}^\eta  + \left\| {{u_2}\left( t \right)} \right\|_{{H^s}}^\eta } \right)}\\{\rm \;\;\;\;\;\;\;\;\;\;\;\;}\le \frac{1}{2}\left\| {{\Lambda ^s}\xi \left( t \right)} \right\|_{{L^2}}^2 + C\left\| {\xi\left( t \right)} \right\|_{{L^2}}^2\left( {1 + \left\| {{u_1}\left( t \right)} \right\|_{{H^s}}^2 + \left\| {{u_2}\left( t \right)} \right\|_{{H^s}}^2} \right)
\end{array}
\end{equation} 
Based on the above discussion we know that $\left\| {{\Lambda ^s}\xi \left( t \right)} \right\|_{{L^2}}^2$ is also bounded, then following the results of \eqref{GrindEQ__5_18_} and \eqref{GrindEQ__5_19_}, yields
\[D_{t}^{\beta } \left\| \xi\left(t\right)\right\| _{L^{2} }^{2} \le K\left(1+\left\| u_{1} \left(t\right)\right\| _{H^{s} }^{2} +\left\| u_{2} \left(t\right)\right\| _{H^{s} }^{2} \right)\left\| \xi\left(t\right)\right\| _{L^{2} }^{2} ,\] 
where $K$ stands for an insignificant constant. Because two different solutions of this equation have the same initial value, namely  $\left\| \xi\left(0\right)\right\| _{L^{2} } =\left\| u_{1} \left(0\right)-u_{2} \left(0\right)\right\| _{L^{2} } =0$. As an application of Gronwall inequality, see Lemma \ref{lemma2.11}, it holds $\left\| \xi\left(t\right)\right\| _{L^{2} } \equiv 0$, which implies the uniqueness of the solution.
\qed

\begin{remark}
We scale and deform fractional inequalities by comparison principle and Property \ref{property2.1}, the similar method to that in \cite{flandoli2021delayed} is used to prove that the lemma holds. From Lemma \ref{lemma5.1}, we can derive the following theorem to hold.
\end{remark}

To prove the main conclusions of this paper, we present several lemmas and corollaries

\begin{lemma}\label{lemma5.2} Let spaces ${\rm{{\cal K}}}$ be $ C\left( {\left[ {0,T} \right];{L^2}} \right) \cap {C^\kappa }\left( {\left[ {0,T} \right];{H^{ - \alpha }}} \right) \cap {L^2}\left( {0,T;{H^s}} \right)$ and ${\rm{{\cal X}}} $ be ${L^p}\left( {0,T;{L^2}} \right) \cap C\left( {\left[ {0,T} \right];{H^{ - \delta }}} \right) \cap {L^2}\left( {0,T;{H^{s - \delta }}} \right)$  \footnote{The definition of the two topologies in this paper are similar to those of Lemma 3.3 in \cite{flandoli2021delayed}.}, where $\alpha$, $\delta$, $s$, $\kappa  >0$, $p < \infty $, then ${\rm {\mathcal K}}\subset\subset {\rm {\mathcal X}}$ (i.e. compact embedding). Moreover, for initial data ${u_M}\left( 0 \right) = {\Pi _M}{u_0}$, the solution of the time fractional order stochastic equation with cut-off function
\[\begin{array}{l}
D_t^\beta {u_M}\left( t \right) = \left[ { - {\Lambda ^{2s}}{u_M}\left( t \right) + b\Delta {u_M}\left( t \right) + {L_S}\left( {{{\left\| {{u_M}\left( t \right)} \right\|}_{{H^{ - \gamma }}}}} \right){\Pi _M}\zeta \left( {{u_M}\left( t \right)} \right)} \right]dt\\
\quad \quad \quad \quad  - A\sum\limits_m {\sum\limits_{j = 1}^{d - 1} {{\theta _m}{\Pi _M}\left( {{\sigma _{m,j}}\cdot\nabla {u_M}\left( t \right)} \right)dW_t^{m,j},} } 
\end{array}\]
satisfy the inequality
\begin{equation}\label{5.20}
\mathbb E\left( {{{\left\| {{u_M}} \right\|}_{{L^2}{H^s}}} + {{\left\| {{u_M}} \right\|}_{{L^\infty }{L^2}}} + {{\left\| {{u_M}} \right\|}_{{C^\kappa }{H^{ - \alpha }}}}} \right) < \infty. 
\end{equation}
\end{lemma}

\noindent $\mathbf{Proof.}$  
Take a sequence ${\left\{ {{h_n}} \right\}_n} \subset {\mathcal K}$ that is bounded. In $ {C^\kappa }\left( {\left[ {0,T} \right];{H^{ - \alpha }}} \right)$, we can extract a subsequence ${\left\{ {{h_{{n_k}}}} \right\}_{{n_k}}}$ from ${\left\{ {{h_n}} \right\}_n} $ that converges to $h$, the conclusion follows from the Ascoli--Arzela theorem \cite{yamazaki2021ascoli}.
By using the inequality \cite{simon1986compact}
\begin{equation}\label{5.1}
{\left\| {{h_{{n_k}}} - h} \right\|_{{H^{ - \delta }}}} \le \left\| {{h_{{n_k}}} - h} \right\|_{{L^2}}^{1 - \mu }\left\| {{h_{{n_k}}} - h} \right\|_{{H^{ - \alpha }}}^\mu 
\end{equation}
where $0 < \mu  < 1$. In space $C\left( {\left[ {0,T} \right];{L^2}} \right)$, it holds
\begin{equation}\label{3}
{\left\| {{h_{{n_k}}}\left( t \right) - h\left( t \right)} \right\|_{{L^2}}} \le C,
\end{equation}
and $C$ is a constant independent of the parameters. Meanwhile, we have 
\begin{equation}\label{4}
{\left\| {{h_{{n_k}}}\left( t \right) - h\left( t \right)} \right\|_{{H^{ - \alpha }}}} \to 0.
\end{equation}
From \eqref{5.1}-\eqref{4}, we can easily deduce that the sequence ${\left\{ {{h_{{n_k}}}} \right\}_{{n_k}}}$ converge to $h$ in space $C\left( {\left[ {0,T} \right];{H^{ - \delta }}} \right)$ for any $\delta  > 0$.  Recall the interpolation inequality, then
\[\begin{array}{l}
\int_0^T {\left\| {{h_{{n_k}}}\left( t \right) - h\left( t \right)} \right\|_{{L^2}}^pdt \le } \int_0^T {\left\| {{h_{{n_k}}}\left( t \right) - h\left( t \right)} \right\|_{{H^\alpha }}^{p\frac{\beta }{{\alpha  + \beta }}}\left\| {{h_{{n_k}}}\left( t \right) - h\left( t \right)} \right\|_{{H^{ - \beta }}}^{p\frac{\alpha }{{\alpha  + \beta }}}dt} \\{\rm \;\;\;\;\;\;\;\;\;\;\;\;\;\;\;\;\;\;\;\;\;\;\;\;\;\;\;\;\;\;\;\;\;\;\;\;\;\;}
 \le \mathop {\sup }\limits_{t \in \left[ {0,T} \right]} \left\| {{h_{{n_k}}}\left( t \right) - h\left( t \right)} \right\|_{{H^{ - \beta }}}^{p\frac{\alpha }{{\alpha  + \beta }}}\int_0^T {\left\| {{h_{{n_k}}}\left( t \right) - h\left( t \right)} \right\|_{{H^\alpha }}^{p\frac{\beta }{{\alpha  + \beta }}}dt},
\end{array}\]
which implies convergence of ${h_{{n_k}}}\left( t \right)$ in the space  ${L^p}\left( {0,T;{L^2}} \right)$. Applying the inequality \eqref{5.1} again after taking the parameter $p = 2$, we state
\[\int_0^T {\left\| {{h_{{n_k}}}\left( t \right) - h\left( t \right)} \right\|_{{H^{s - \delta }}}^2dt}  \le \int_0^T {\left\| {{h_{{n_k}}}\left( t \right) - h\left( t \right)} \right\|_{{H^s}}^{2\left( {1 - \theta } \right)}\left\| {{h_{{n_k}}}\left( t \right) - h\left( t \right)} \right\|_{{L^2}}^{2\theta }dt}. \]
Similarly, since ${h_{{n_k}}}\left( t \right)$ converge to $ h\left( t \right)$ in ${L^2}\left( {0,T;{L^2}} \right)$ and the boundness of ${L^2}\left( {0,T;{H^s}} \right)$, the sequence converge in space ${L^2}\left( {0,T;{H^s}} \right)$. 

Next, we use the notation in Remark \ref{remark5.1} and verify the second conclusion. If $u \in {H^s}$, then $\Delta u \in {H^{s - 2}} \subset {H^{ - s}}$, according to the triangle inequality,
\[\begin{array}{l}
\frac{1}{{{\Gamma ^2}\left( \beta  \right)}}\int_0^T {\left\| {{{\left( {T - \tau } \right)}^{\beta  - 1}}\left[ { - {{\left( { - \Delta } \right)}^s}u\left( \tau  \right) + b\Delta u\left( \tau  \right) + {L_S}\left( {{{\left\| {u\left( \tau  \right)} \right\|}_{{H^{ - \gamma }}}}} \right)\zeta \left( {u\left( \tau  \right)} \right)} \right]} \right\|_{{H^{ - s}}}^2} d\tau \\{\rm \;\;\;\;\;\;\;\;\;\;\;\;\;}
 \le \frac{1}{{{\Gamma ^2}\left( \beta  \right)}}\int_0^T {{{\left( {T - \tau } \right)}^{2\beta  - 2}}\left\| { - {{\left( { - \Delta } \right)}^s}u\left( \tau  \right) + b\Delta u\left( \tau  \right)} \right\|_{{H^{ - s}}}^2} d\tau \\{\rm \;\;\;\;\;\;\;\;\;\;\;\;\;}
 + \frac{1}{{{\Gamma ^2}\left( \beta  \right)}}\int_0^T {{{\left( {T - \tau } \right)}^{2\beta  - 2}}\left\| {{L_S}\left( {{{\left\| {u\left( \tau  \right)} \right\|}_{{H^{ - \gamma }}}}} \right)\zeta \left( {u\left( \tau  \right)} \right)} \right\|_{{H^{ - s}}}^2} d\tau \\{\rm \;\;\;\;\;\;\;\;\;\;\;\;\;}
 \lesssim \frac{1}{{{\Gamma ^2}\left( \beta  \right)}}\int_0^T {{{\left( {T- \tau } \right)}^{2\beta  - 2}}\left( {1 + {b^2}} \right)\left\| {u\left( \tau  \right)} \right\|_{{H^s}}^2} d\tau  + \frac{1}{{{\Gamma ^2}\left( \beta  \right)}}\int_0^T {{{\left( {T - \tau } \right)}^{2\beta  - 2}}\left\| {\zeta \left( {u\left( \tau  \right)} \right)} \right\|_{{H^{ - s}}}^2} d\tau,
\end{array}\]
by using the Hypothesis \ref{hypothesis3.1} and H{\"o}lder's inequality, it holds
\[\frac{1}{{{\Gamma ^2}\left( \beta  \right)}}\int_0^T {{{\left( {T - \tau } \right)}^{2\beta  - 2}}\left( {1 + {b^2}} \right)\left\| {u\left( \tau  \right)} \right\|_{{H^s}}^2} d\tau  \le C_1\left\| u \right\|_{{L^2}{H^s}}^2.\]
Similarly, we can easily show that
\[\frac{1}{{{\Gamma ^2}\left( \beta  \right)}}\int_0^T {{{\left( {T - \tau } \right)}^{2\beta  - 2}}\left\| {\zeta \left( {u\left( \tau  \right)} \right)} \right\|_{{H^{ - s}}}^2} d\tau  \le C_2\left( {1 + \left\| u \right\|_{{L^2}{H^s}}^2} \right)\left( {1 + \left\| u \right\|_{{L^\infty }{L^2}}^{2{a_1}}} \right),\]
where $C_1,C_2$ are both constants. Therefore, we deduce
\[{\left\| {G\left( t \right)} \right\|_{{C^{1/2}}{H^{ - s}}}} \le {\left\| {G\left( t \right)} \right\|_{{W^{1,2}}{H^{ - s}}}} \mathbin{\lesssim} \left( {1 + \left\| u \right\|_{{L^\infty }{L^2}}^{{a_1}}} \right)\left( {1 + {{\left\| u \right\|}_{{L^2}{H^s}}}} \right).\]
For $M\left( t \right)$, we have
\[\begin{array}{l}
{\left[ {\left\langle {M\left( t \right),{e^{2\pi im\cdot x}}} \right\rangle } \right]\left( t \right) = \frac{{{A^2}}}{{\Gamma {{\left( \beta  \right)}^2}}}\left[ {\sum\limits_m {\sum\limits_{j = 1}^{d - 1} {{\theta _m}\int_0^\cdot {\left\langle {{{\left( {t - r} \right)}^{\beta  - 1}}{u_M}\left( r \right),{\sigma _{k,i}}\cdot{\Pi _M}\nabla {e^{2\pi im\cdot x}}} \right\rangle } dW_r^{m,j}} } } \right]\left( t \right)}\\{\rm \;\;\;\;\;\;\;\;\;\;\;\;\;}{ \le \frac{{2{A^2}}}{{\Gamma {{\left( \beta  \right)}^2}}}\sum\limits_m {\sum\limits_{j = 1}^{d - 1} {\theta _m^2\int_0^t {{{\left( {t - r} \right)}^{2\beta  - 2}}{{\left\langle {{u_M}\left( r \right),{\sigma _{m,j}}\cdot{\Pi _M}\nabla {e^{2\pi im\cdot x}}} \right\rangle }^2}} dr} } }\\{\rm \;\;\;\;\;\;\;\;\;\;\;\;\;}{ \le \frac{{2{A^2}\left\| \theta  \right\|_{{\ell ^\infty }}^2}}{{\Gamma {{\left( \beta  \right)}^2}}}\int_0^t {{{\left( {t - r} \right)}^{2\beta  - 2}}\left| {\left| {{u_M}\left( r \right){\Pi _M}\nabla {e^{2\pi im\cdot x}}} \right|} \right|_{{L^2}}^2dr} }\\{\rm \;\;\;\;\;\;\;\;\;\;\;\;\;}{ \le \frac{{2{A^2}\left\| \theta  \right\|_{{\ell ^\infty }}^2}}{{\Gamma {{\left( \beta  \right)}^2}}}\left| {\left| {{\Pi _M}\nabla {e^{2\pi im\cdot x}}} \right|} \right|_{{L^\infty }}^2\int_0^t {{{\left( {t - r} \right)}^{2\beta  - 2}}\left| {\left| {{u_M}\left( r \right)} \right|} \right|_{{L^2}}^2dr} }.
\end{array}\]

Since any order partial derivative of ${e^{2\pi im \cdot x}}$  exists,
$\nabla {e^{2\pi im \cdot x}} = 2\pi mi \cdot {e^{2\pi im \cdot x}}$, 
$\Delta {e^{2\pi im \cdot x}} =  - 4{\pi ^2}{m^2} \cdot {e^{2\pi im \cdot x}}$, ${\sigma _{m,j}} = {q_{m,j}}{e^{2\pi im \cdot x}}$. We shall use the Burkholder-Davis-Gundy inequality \cite{ren2008burkholder} and H{\"o}lder's inequality to estimate the term $M(t)$, it holds
\[\mathbb E\left[ {\left\| {M\left( t \right)} \right\|_{{H^{ - \alpha }}}^2} \right] \lesssim \sum\limits_m {{{\left( {1 + {{\left| m \right|}^2}} \right)}^{ - \alpha }}E\left[ {{{\left| {\left\langle {{M_M}\left( t \right),{e^{2\pi im\cdot x}}} \right\rangle } \right|}^2}} \right] \lesssim \left\| \theta  \right\|_{{\ell ^\infty }}^2{t^{2\beta  - 1}}} \]
By using the Lemma \ref{lemma5.1}, it is obvious that the inequality \eqref{5.2} holds. \qed

\begin{remark}
This lemma is essential to prove the tightness of the distribution of solutions to the Eq.\eqref{GrindEQ__5_71_}. We use the similar method in \cite{flandoli2021delayed} to prove the embedding relation of two topologies. Based on this lemma we can get the following lemma to hold.
\end{remark}

\begin{corollary}\label{corollary1} With ${\rho _M}$ denoting the law of ${{u_M}}$, the family ${\left\{ {{\rho _M}} \right\}_M}$ is tight on 
\[{\rm{{\cal X}}} = {L^p}\left( {0,T;{L^2}} \right) \cap C\left( {\left[ {0,T} \right];{H^{ - \delta }}} \right) \cap {L^2}\left( {0,T;{H^{s - \delta }}} \right).\]
\end{corollary}
\noindent $\mathbf{Proof.}$ Recall the conclusion of Lemma \ref{lemma5.2} that $K \subset  \subset X$, then \[{\left\| {\; \cdot \;} \right\|_{\mathcal{X}}} \lesssim {\left\| {\; \cdot \;} \right\|_\mathcal{K}},\] it is obvious from the Prokhorov theorem (see \cite{chin2019new}) that the conclusion holds.  \qed

\begin{proposition}\label{pro2} Take the sequence $\left\{\theta ^{S} \right\}_{S\ge 1} \subset \ell ^{2} $ that satisfy the limit ${{{{\left\| {{\theta ^S}} \right\|}_{{\ell ^\infty }}}} \mathord{\left/
 {\vphantom {{{{\left\| {{\theta ^S}} \right\|}_{{\ell ^\infty }}}} {{{\left\| {{\theta ^S}} \right\|}_{{\ell ^2}}}}}} \right.\kern-\nulldelimiterspace} {{{\left\| {{\theta ^S}} \right\|}_{{\ell ^2}}}}} \to 0$ as $S \to \infty$. The bounded, convex and closed set ${\left\{ {u_0^S} \right\}_S} \in {L^2}\left( {{\mathbb T^d}} \right)$ satisfies the condition ${u_0^S}\stackrel{W}{\longrightarrow}{u_0}$. On ${L^p}\left( {0,T;{L^2}} \right) \cap C\left( {\left[ {0,T} \right];{H^{ - \delta }}} \right) \cap {L^2}\left( {0,T;{H^{s - \delta }}} \right)$, the unique solution of fractional stochastic differential equation 
\begin{equation}\label{5.2}
\left\{ {\begin{array}{*{20}{l}}
{D_t^\beta u = \left[ { - {\Lambda ^{2s}}u\left( \tau  \right) + b\Delta u\left( \tau  \right) + {L_S}\left( {{{\left\| u \right\|}_{{H^{ - \gamma }}}}} \right)\zeta \left( {u\left( \tau  \right)} \right)} \right]dt + A\sum\limits_m {\sum\limits_{j = 1}^{d - 1} {\theta _m^S{\sigma _{m,j}}\cdot\nabla u\left( t \right)dW_t^{m,j}} } }\\
{{u_0} = u_0^S}
\end{array}} \right.
\end{equation}
converge to the unique solution of the time fractional deterministic equation 
\begin{equation}\label{5.3}
D_t^\beta u = \left[ { - {\Lambda ^{2s}}u\left( \tau  \right) + b\Delta u\left( \tau  \right) + {L_S}\left( {{{\left\| u \right\|}_{{H^{ - \gamma }}}}} \right)\zeta \left( {u\left( \tau  \right)} \right)} \right]dt,
\end{equation}
in probability with initial data ${u_0}$.
\end{proposition}

\noindent $\mathbf{Proof}$ The symbol ${\left\{ {{\rho ^S}} \right\}_S}$ is used to represent the law of ${\left\{ {{u^S}} \right\}_S}$. It follows from the Prokhorov theorem and above results that exists a subsequence ${\left\{ {{\rho ^{{S_i}}}} \right\}_{{S_i}}} \subset {\left\{ {{\rho ^S}} \right\}_S}$, and holds ${\rho ^{{S_i}}}\stackrel{W}{\longrightarrow}\rho$, while $\rho$ is a probability  measure. Then, Skorokhod's representation theorem (see \cite{cortissoz2007skorokhod}) tells us, on another probability space $\left( {\mathord{\buildrel{\lower3pt\hbox{$\scriptscriptstyle\frown$}} 
\over \Omega } ,\mathord{\buildrel{\lower3pt\hbox{$\scriptscriptstyle\frown$}} 
\over {\mathcal{F}}} ,\mathord{\buildrel{\lower3pt\hbox{$\scriptscriptstyle\frown$}} 
\over {\mathbb P}}} \right)$, it is possible to find the random variables ${\left\{ {{{\mathord{\buildrel{\lower3pt\hbox{$\scriptscriptstyle\frown$}} 
\over u} }^{{S_i}}}} \right\}_{{S_i}}}$ and $\mathord{\buildrel{\lower3pt\hbox{$\scriptscriptstyle\frown$}} 
\over u} $, the former have the same distribution as ${\rho ^{{S_i}}}$ and the latter  have the same distribution as $\rho $, such that ${\mathord{\buildrel{\lower3pt\hbox{$\scriptscriptstyle\frown$}} 
\over u} ^{{S_i}}}\stackrel{a.s.}{\longrightarrow}\mathord{\buildrel{\lower3pt\hbox{$\scriptscriptstyle\frown$}} 
\over u}$.

Next, we consider the new complex motions, which are defined on the new probability space $\left( {\mathord{\buildrel{\lower3pt\hbox{$\scriptscriptstyle\frown$}} 
\over \Omega } ,\mathord{\buildrel{\lower3pt\hbox{$\scriptscriptstyle\frown$}} 
\over {\mathcal{F}}} ,\mathord{\buildrel{\lower3pt\hbox{$\scriptscriptstyle\frown$}} \over {\mathbb{P}}} } \right)$, and write as  $\mathord{\buildrel{\lower3pt\hbox{$\scriptscriptstyle\frown$}} 
\over W}  = \left\{ {{{\mathord{\buildrel{\lower3pt\hbox{$\scriptscriptstyle\frown$}} 
\over W} }^{m,j}}:m \in Z_0^d,j = 1, \cdots ,d - 1} \right\}$. It hold 
$\left( {{u^{{S_i}}},W} \right)\ \overset{d}{=} \left( {{{\mathord{\buildrel{\lower3pt\hbox{$\scriptscriptstyle\frown$}} 
\over u} }^{{S_i}}},{{\mathord{\buildrel{\lower3pt\hbox{$\scriptscriptstyle\frown$}} 
\over W} _{{S_i}}}}} \right),$ in addition 
$\mathord{\buildrel{\lower3pt\hbox{$\scriptscriptstyle\frown$}} 
\over W} _{{S_i}}^{m,j}\stackrel{a.s.}{\longrightarrow}\mathord{\buildrel{\lower3pt\hbox{$\scriptscriptstyle\frown$}} 
\over W} ^{m,j}.$ Note that $\mathord{\buildrel{\lower3pt\hbox{$\scriptscriptstyle\frown$}} 
\over W} _{{S_i}}^{m,j}$ represents the subsequence of ${\left\{ {{{\mathord{\buildrel{\lower3pt\hbox{$\scriptscriptstyle\frown$}} 
\over W} }^{m,j}}:m \in Z_0^d,j = 1, \cdots ,d} \right\}}$, other conditions are the same as the family ${\left\{ {W_t^{m,j}} \right\}_{m,j}}$.

We shall prove that $\mathord{\buildrel{\lower3pt\hbox{$\scriptscriptstyle\frown$}} 
\over u}$ is the solution to deterministic equation \eqref{5.3}.
Obviously,
\[\begin{array}{l}
\left\langle {{{\mathord{\buildrel{\lower3pt\hbox{$\scriptscriptstyle\frown$}} 
\over u} }^{{S_i}}}\left( t \right),\varphi } \right\rangle  = \left\langle {u_0^{{S_i}},\varphi } \right\rangle  - \frac{1}{{\Gamma \left( \beta  \right)}}\int_0^t {\left[ {\left\langle {{{\left( {t - \tau } \right)}^{\beta  - 1}}{\Delta ^s}{{\mathord{\buildrel{\lower3pt\hbox{$\scriptscriptstyle\frown$}} 
\over u} }^{{S_i}}}\left( \tau  \right),\varphi } \right\rangle  + b\left\langle {{{\left( {t - \tau } \right)}^{\beta  - 1}}\Delta {{\mathord{\buildrel{\lower3pt\hbox{$\scriptscriptstyle\frown$}} 
\over u} }^{{S_i}}}\left( \tau  \right),\varphi } \right\rangle } \right]d\tau } \\
{\rm{   \;\;\;\;\;\;\;\;\;\;\;\;\;\;\;\;\;\;        }} + \frac{1}{{\Gamma \left( \beta  \right)}}\int_0^t {{L_S}\left( {{{\left\| {{{\mathord{\buildrel{\lower3pt\hbox{$\scriptscriptstyle\frown$}} 
\over u} }^{{S_i}}}} \right\|}_{{H^{ - \gamma }}}}} \right)\left\langle {{{\left( {t - \tau } \right)}^{\beta  - 1}}\zeta \left( {{{\mathord{\buildrel{\lower3pt\hbox{$\scriptscriptstyle\frown$}} 
\over u} }^{{S_i}}}\left( \tau  \right)} \right),\varphi } \right\rangle } d\tau \\
{\rm{   \;\;\;\;\;\;\;\;\;\;\;\;\;\;\;\;\;\;          }} - \frac{A}{{\Gamma \left( \beta  \right)}}\sum\limits_m {\sum\limits_{j = 1}^{d - 1} {\theta _m^{{S_i}}\int_0^t {\left\langle {{{\left( {t - \tau } \right)}^{\beta  - 1}}{{\mathord{\buildrel{\lower3pt\hbox{$\scriptscriptstyle\frown$}} 
\over u} }^{{S_i}}}\left( \tau  \right),{\sigma _{m,j}} \cdot \nabla \varphi } \right\rangle } dW_\tau ^{m,j}} } 
\end{array}\]
for test function $\varphi  \in {C^\infty }\left( {{\mathbb T^d},{\mathbb R^d}} \right)$, which is divergence free.  For the last part of the above inequality, we have \footnote{Drawing from the method of Proposition 4.2 in \cite{flandoli2021high} for estimating the martingale part, the conclusion of this paper can be seen as a special case when the parameter $\beta$ is not equal to 1.} 
\[\begin{array}{l}
\mathord{\buildrel{\lower3pt\hbox{$\scriptscriptstyle\frown$}} 
\over {\mathbb E}} {\left( {{M^{{S_i}}}\left( t \right)} \right)^2} = \mathord{\buildrel{\lower3pt\hbox{$\scriptscriptstyle\frown$}} 
\over {\mathbb E}} {\left( {\frac{A}{{\Gamma \left( \beta  \right)}}\sum\limits_m {\sum\limits_{j = 1}^{d - 1} {\theta _m^{{S_i}}\int_0^t {\left\langle {{{\left( {t - \tau } \right)}^{\beta  - 1}}{{\mathord{\buildrel{\lower3pt\hbox{$\scriptscriptstyle\frown$}} 
\over u} }^{{S_i}}}\left( \tau  \right),{\sigma _{m,j}} \cdot \nabla \varphi } \right\rangle } dW_\tau ^{m,j}} } } \right)^2}\\{\rm \;\;\;\;\;\;\;\;\;\;\;\;\;\;\;\;\;\;\;\;}
 = 2\frac{{{A^2}}}{{\Gamma {{\left( \beta  \right)}^2}}}\sum\limits_m {\sum\limits_{j = 1}^{d - 1} {{{\left( {\theta _m^{{S_i}}} \right)}^2}\mathord{\buildrel{\lower3pt\hbox{$\scriptscriptstyle\frown$}} 
\over {\mathbb E}} \int_0^t {{{\left| {\left\langle {{{\left( {t - \tau } \right)}^{\beta  - 1}}{{\mathord{\buildrel{\lower3pt\hbox{$\scriptscriptstyle\frown$}} 
\over u} }^{{S_i}}}\left( \tau  \right),{\sigma _{m,j}} \cdot \nabla \varphi } \right\rangle } \right|}^2}} d\tau } } \\{\rm \;\;\;\;\;\;\;\;\;\;\;\;\;\;\;\;\;\;\;\;}
 \lesssim \frac{{{A^2}\left\| \theta  \right\|_{{\ell ^\infty}}^2 }}{{\Gamma {{\left( \beta  \right)}^2}}}\mathord{\buildrel{\lower3pt\hbox{$\scriptscriptstyle\frown$}} 
\over {\mathbb E}} \int_0^t {\sum\limits_m {\sum\limits_{j = 1}^{d - 1} {{{\left| {\left\langle {{{\left( {t - \tau } \right)}^{\beta  - 1}}{{\mathord{\buildrel{\lower3pt\hbox{$\scriptscriptstyle\frown$}} 
\over u} }^{{S_i}}}\left( \tau  \right),{\sigma _{m,j}} \cdot \nabla \varphi } \right\rangle } \right|}^2}} } d\tau },
\end{array}\]
where ${\mathord{\buildrel{\lower3pt\hbox{$\scriptscriptstyle\frown$}} 
\over {\mathbb E}} }$ represent the expectation on the new probability space $\left( {\mathord{\buildrel{\lower3pt\hbox{$\scriptscriptstyle\frown$}} 
\over \Omega } ,\mathord{\buildrel{\lower3pt\hbox{$\scriptscriptstyle\frown$}} 
\over {\mathcal{F}}} ,\mathord{\buildrel{\lower3pt\hbox{$\scriptscriptstyle\frown$}} 
\over {\mathbb P}} } \right)$. By recalling the definition of
${\left\{ {{\sigma _{m,j}}} \right\}_{m,j}} \in {L^2}\left( {{\mathbb T^d},{\mathbb R^d}} \right)$, we show that
\[\begin{array}{l}
{\sum\limits_m {\sum\limits_{j = 1}^{d - 1} {\left| {\left\langle {{{\left( {t - \tau } \right)}^{\beta  - 1}}{{\mathord{\buildrel{\lower3pt\hbox{$\scriptscriptstyle\frown$}} 
\over u} }^{{S_i}}}\left( \tau  \right),{\sigma _{m,j}} \cdot \nabla \varphi } \right\rangle } \right|} } ^2} = {\sum\limits_m {\sum\limits_{j = 1}^{d - 1} {\left| {\left\langle {\nabla \varphi  \cdot {{\left( {t - \tau } \right)}^{\beta  - 1}}{{\mathord{\buildrel{\lower3pt\hbox{$\scriptscriptstyle\frown$}} 
\over u} }^{{S_i}}}\left( \tau  \right),{\sigma _{m,j}}} \right\rangle } \right|} } ^2}\\{\rm \;\;\;\;\;\;\;\;\;\;\;\;}
 \le \left\| {\nabla \varphi  \cdot {{\left( {t - \tau } \right)}^{\beta  - 1}}{{\mathord{\buildrel{\lower3pt\hbox{$\scriptscriptstyle\frown$}} 
\over u} }^{{S_i}}}\left( \tau  \right)} \right\|_{{L^2}}^2 \le \left\| {\nabla \varphi } \right\|_{{L^\infty }}^2\left\| {{{\left( {t - \tau } \right)}^{\beta  - 1}}{{\mathord{\buildrel{\lower3pt\hbox{$\scriptscriptstyle\frown$}} 
\over u} }^{{S_i}}}\left( \tau  \right)} \right\|_{{L^2}}^2.
\end{array}\]
According to the value of the parameter $A$ and the conditions satisfied by the special sequence ${\left\{ {{\theta ^{{S_i}}}} \right\}_{i \ge 1}}$, it is obvious that when $i$ tends to infinity,
\[\frac{{{A^2}\left\| \theta  \right\|_{{\ell ^\infty }}^2}}{{\Gamma {{\left( \beta  \right)}^2}}} = \frac{{d \cdot b}}{{d - 1}}\frac{{\left\| \theta  \right\|_{{\ell ^\infty }}^2}}{{\left\| \theta  \right\|_{{\ell ^2}}^2}} \to 0.\]
From the discussion above, the boundedness of  ${\left\| {{{\mathord{\buildrel{\lower3pt\hbox{$\scriptscriptstyle\frown$}} 
\over u} }^{{S_i}}}} \right\|_{{L^2}}}$ is known, then combining with the H{\"o}lder's inequality, we deduce the expectation
$\mathord{\buildrel{\lower3pt\hbox{$\scriptscriptstyle\frown$}} 
\over {\mathbb E}} {\left( {{M^{{S_i}}}\left( t \right)} \right)^2} $ tends to $0$ as $i$ tends to infinity. Recall the fact ${\mathord{\buildrel{\lower3pt\hbox{$\scriptscriptstyle\frown$}} 
\over u} ^{{S_i}}}\stackrel{a.s.}{\longrightarrow}\mathord{\buildrel{\lower3pt\hbox{$\scriptscriptstyle\frown$}} 
\over u}$ and $u_0^S\stackrel{W}{\longrightarrow}{u_0}$, we can conclude the following identity holds,
\[\begin{array}{l}
\left\langle {\mathord{\buildrel{\lower3pt\hbox{$\scriptscriptstyle\frown$}} 
\over u} \left( t \right),\varphi } \right\rangle  = \left\langle {{u_0},\varphi } \right\rangle  - \frac{1}{{\Gamma \left( \beta  \right)}}\int_0^t {{{\left( {t - \tau } \right)}^{\beta  - 1}}\left[ {\left\langle {{\Lambda ^s}\mathord{\buildrel{\lower3pt\hbox{$\scriptscriptstyle\frown$}} 
\over u} \left( \tau  \right),{\Lambda ^s}\varphi } \right\rangle  + b\left\langle {\nabla \mathord{\buildrel{\lower3pt\hbox{$\scriptscriptstyle\frown$}} 
\over u} \left( \tau  \right),\nabla \varphi } \right\rangle } \right]d\tau } \\
{\rm{ \;\;\;\;\;\;\;\;\;\;\;\;\;          }} + \frac{1}{{\Gamma \left( \beta  \right)}}\int_0^t {{{\left( {t - \tau } \right)}^{\beta  - 1}}{L_S}\left( {{{\left\| {\mathord{\buildrel{\lower3pt\hbox{$\scriptscriptstyle\frown$}} 
\over u} } \right\|}_{{H^{ - \gamma }}}}} \right)\left\langle {\zeta \left( {\mathord{\buildrel{\lower3pt\hbox{$\scriptscriptstyle\frown$}} 
\over u} \left( \tau  \right)} \right),\varphi } \right\rangle } d\tau.
\end{array}\]
Namely, $\mathord{\buildrel{\lower3pt\hbox{$\scriptscriptstyle\frown$}} \over u}$ is the solution to deterministic equation \eqref{5.3}. 
Overall results of the above, the family ${\left\{ {{\rho ^{{S_i}}}} \right\}_{{S_i}}}$ converge weakly to ${\delta _u}$, which represents the Dirac measure, also known as the point mass at $u$, note that $u$ represents the unique solution of deterministic equation. From what we know from corollary \ref{corollary1}, the ${\left\{ {{\rho ^S}} \right\}_S}$  is tight on $\mathcal{X}$, from which we claim that ${\left\{ {{\rho ^S}} \right\}_S}\stackrel{W}{\longrightarrow}{\delta _u}$.  As the sequence of distribution functions  converges weakly, the sequence of random variables ${\left\{ {{u^S}} \right\}_S}$ converges in distribution to $u$, so moreover, the sequence also converges in probability.   \qed

\begin{remark}
Inspired by \cite{flandoli2021delayed} and \cite{flandoli2021high}, we prove the convergence of the solution to the fractional stochastic equation \eqref{5.2},  which can be seen as an extension of the conclusion about the fractional-order equation and of vital importance for the following proof of the main theorem in this paper.
\end{remark}

\section{Proof of main results}\label{proof}
In this section, we shall give proof of the two most important conclusions of this paper.

\noindent $\mathbf{Proof{\rm \; \;} of{\rm \; \;} Theorem{\rm \; \;} 2.1.}$ Recall the condition (\romannumeral4) in the Hypothesis \ref{hypothesis3.1}, for the bounded ${u_0}$, it is possible to find constant $S$ and $b$ that 
\[{\left\| {u_t}\left( {{u_0}} \right) \right\|_{{L^\infty }\left( {0,T;{L^2}} \right)}} \le S - 1,\]
By Proposition \ref{pro2}, it is clear that 
\[\mathbb P\left( {{{\left\| {u_t^S\left( {{u_0},\theta } \right) - {u_t}\left( {{u_0}} \right)} \right\|}_{C\left( {\left[ {0,T} \right];{H^{ - \gamma }}} \right)}} \le \varepsilon } \right) \ge 1 - \varepsilon ,\]
for any $\varepsilon  \in \left( {0,1} \right)$ and $T > 0$.
Apply the properties of the norm, yield
\[{\left\| {u_t^S\left( {{u_0},{\theta ^S}} \right)} \right\|_{C\left( {\left[ {0,T} \right];{H^{ - \gamma }}} \right)}} - {\left\| {{u_t}\left( {{u_0}} \right)} \right\|_{C\left( {\left[ {0,T} \right];{H^{ - \gamma }}} \right)}} \le {\left\| {u_t^S\left( {{u_0},{\theta ^S}} \right) - {u_t}\left( {{u_0}} \right)} \right\|_{C\left( {\left[ {0,T} \right];{H^{ - \gamma }}} \right)}},\]
we can deduce that ${\left\| {u_t^S\left( {{u_0},{\theta ^S}} \right)} \right\|_{C\left( {\left[ {0,T} \right];{H^{ - \gamma }}} \right)}} < S,$ which implies that
\[\mathbb P\left( {{{\left\| {u_t^S\left( {{u_0},{\theta ^S}} \right)} \right\|}_{C\left( {\left[ {0,T} \right];{H^{ - \gamma }}} \right)}} < S} \right) \ge P\left( {{{\left\| {u_t^S\left( {{u_0},{\theta ^S}} \right) - {u_t}\left( {{u_0}} \right)} \right\|}_{C\left( {\left[ {0,T} \right];{H^{ - \gamma }}} \right)}} \le \varepsilon } \right) \ge 1 - \varepsilon. \]
in others word, 
\[\mathbb P\left( {{L_S}\left( {{{\left\| {u_t^S\left( {{u_0},{\theta ^S}} \right)} \right\|}_{C\left( {\left[ {0,T} \right];{H^{ - \gamma }}} \right)}}} \right) = 1} \right) \ge 1 - \varepsilon ,\]
it also means that the  Eq. \eqref{5.3} reduces to the deterministic equation without cut-off. The above result implies that the life span of solutions to \eqref{5.2} with initial data $u_{0} $ is greater than $T$.    \qed

\begin{remark}
Compared with \cite{flandoli2021high}, the proof of Theorem \ref{theorem3.1} in this paper mainly employs the conclusion of Proposition \ref{pro2}, the triangle inequality and the inclusion relation of events to prove the boundedness of the norm of the solution to Eq. \eqref{5.3},  and extends the result to the lifetime problem of the solution.
\end{remark}

\noindent $\mathbf{Proof{\rm \; \;} of{\rm \; \;} Theorem{\rm \; \;} 2.2.}$ Since we have the inequality
\[{\left\| {{u_t}\left( {{u_0}} \right)} \right\|_{{L^2}}} \le K{\left\| {{u_0}} \right\|_{{L^2}}}{e^{ - \lambda t}},\]
where parameters $K$ and $\lambda$ can be taken big enough. Selecting the appropriate parameters $K$ and $\lambda$ makes
\[{\left\| {{u_t}\left( {{u_0}} \right)} \right\|_{{L^2}\left( {T - 1,T;{L^2}} \right)}} \le {\left[ {\int_{T-1}^T {{{\left( {K{{\left\| {{u_0}} \right\|}_{{L^2}}}{e^{ - \lambda t}}} \right)}^2}dt} } \right]^{{1 \mathord{\left/
 {\vphantom {1 2}} \right.
 \kern-\nulldelimiterspace} 2}}} \le \frac{{{K^2}\left\| {{u_0}} \right\|_{{L^2}}^2}}{{2\lambda }}{e^{ - 2\lambda \left( {T - 1} \right)}} \le \frac{\varepsilon }{2}.\]

Together with the conclusion of Proposition \ref{pro2}, it is possible to find $S$ big enough such that the equation  reduce to \eqref{GrindEQ__3_71_}, moreover, it holds
\[\mathbb P\left( {{{\left\| {u_t^S\left( {{u_0},\theta } \right) - {u_t}\left( {{u_0}} \right)} \right\|}_{{L^2}\left( {0,T;{L^2}} \right)}} \le \varepsilon } \right) \ge 1 - \varepsilon ,\]
Let ${\rm Z} = \left\{ {{{\left\| {u_t^S\left( {{u_0},\theta } \right) - {u_t}\left( {{u_0}} \right)} \right\|}_{{L^2}\left( {0,T;{L^2}} \right)}} \le \varepsilon } \right\}$ be a set, namely, the probability of this event occurring is greater than $ 1 - \varepsilon$. It is easily to deduce from triangle inequality that
\[\begin{array}{l}
{\left\| {u_t^S\left( {{u_0},\theta } \right)} \right\|_{{L^2}\left( {T - 1,T;{L^2}} \right)}} \le {\left\| {u_t^S\left( {{u_0},\theta } \right) - {u_t}\left( {{u_0}} \right)} \right\|_{{L^2}\left( {T - 1,T;{L^2}} \right)}} + {\left\| {{u_t}\left( {{u_0}} \right)} \right\|_{{L^2}\left( {T - 1,T;{L^2}} \right)}}\\{\rm \;\;\;\;\;\;\;\;\;\;\;\;\;\;\;\;\;\;\;\;\;\;\;\;\;\;\;\;\;\;\;\;\;\;\;\;\;\;}
 \le \frac{\varepsilon }{2} + \frac{\varepsilon }{2} = \varepsilon,
\end{array}\]
namely we can find $t\left( \upsilon  \right)$, which is taken in the interval  $\left[ {T - 1,T} \right]$ and $\upsilon  \in {\rm Z}$, such that the ${L^2}$-norm of ${u_{t\left( \upsilon  \right)}^S\left( {{u_0},\theta } \right)}$ small enough. Let ${\left\| {u_{t\left( \upsilon  \right)}^S\left( {{u_0},\theta } \right)} \right\|_{{L^2}}}$ be the initial data of Eq.\eqref{GrindEQ__3_71_}, then the conclusion holds. \qed

\section*{Future work}
In this paper, we mainly focus on the fractional stochastic equation
\[\left\{\begin{array}{l} {D_t^\beta u = \left[ { - {{\left( { - \Delta } \right)}^s}u + \zeta \left( u \right)} \right]dt + A\sum\limits_{m \in Z_0^d} {\sum\limits_{j = 1}^{d - 1} {{\theta _m}{\sigma _{m,j}} \circ dW_t^{m,j}} } ,} \\ {u\left|_{t=0} =u_{0} ,\right. } \end{array}\right. \]
considering the uniqueness of its solution and the explosion phenomenon, as well as the influence of noise on its blow-up time. In the subsequent study, we will consider the blow-up problem in equation 
\[\left\{ {\begin{array}{*{20}{l}}
{D_t^{{\beta _1}}u = \left[ { - {{\left( { - \Delta } \right)}^s}u + \zeta \left( u \right)} \right]dt + A\int_0^t {\sum\limits_{m \in Z_0^d} {\sum\limits_{j = 1}^{d - 1} {{\theta _m}{\sigma _{m,j}} \circ D_\tau ^{{\beta _2}}W_\tau ^{m,j}} } } ,}\\
{u\left| {_{t = 0} = {u_0},} \right.}
\end{array}} \right.\]
where ${\beta _1},{\beta _2} \in \left( {0,1} \right)$.
\section*{Conflict of interest}

The authors declare that they have no conflict of interest.

\section*{Acknowledgment}
We would like to thank an anonymous reviewer for providing useful suggestions for improving the presentation of this work.

\begin{appendices}
\section{Definitions and complements} \label{appendixA}
\subsection{Definitions} 

Let $p \le  + \infty $ is a positive integer,  ${L^p}\left( {{{\mathbb T}^d},{{\mathbb R}^d}} \right)$ represents the set of  $p$-th integrable function \cite{ju2004existence}, which is defined on ${\mathbb T}^d$ and real-valued and has norm
   \begin{gather*}
       \left\| f\right\| _{L^{p} } =\left(\int _{\Omega }\left|f\left(x\right)\right|^{p} dx \right)^{\frac{1}{p} } ,{\rm \; \; \; \; \; \; \; \; \; \; }\left\| f\right\| _{L^{\infty } } =ess\mathop{\sup }\limits_{x\in \Omega } \left|f\left(x\right)\right|.
   \end{gather*}
\noindent And the Sobolev space is given in \cite{ju2004existence,galeati2020convergence}, namely
 \begin{gather*}
     H^{s} \left({\mathbb T}^{d} \right)=\left\{f=\sum _{k}f_{k} e_{k}  \left|f_{-k} \right. =\bar{f}_{k} ,\sum _{k}\left(1+\left|k\right|^{2} \right)^{s } \left|f_{k} \right|^{2} <\infty  \right\},
   \end{gather*} 
where $f\in L^{2} \left({\mathbb T}^{d} ;{\mathbb R}\right)$, $\left\{e_{k} \right\}_{k\in {\mathbb Z}^{d} } $ indicates a complete orthonormal system, note that $e_{k} =e^{ik\cdot x} $. Moreover, there is an important property of such Sobolev spaces, ${\left\| f \right\|_{{H^s}}} = {\left\| {{{\left( { - \Delta } \right)}^{{s \mathord{\left/{\vphantom {s 2}} \right.\kern-\nulldelimiterspace} 2}}}f} \right\|_{{L^2}}}$. 

The fractional Sobolev space (\cite{borikhanov2022qualitative})
\[{W^{s,p}}\left( {{\mathbb R^n}} \right): = \left\{ {u \in {L^p}\left( {{\mathbb R^n}} \right),\frac{{\left| {u\left( x \right) - u\left( y \right)} \right|}}{{{{\left| {x - y} \right|}^{\frac{N}{p} + s}}}} \in {L^p}\left( {{R^n} \times {R^n}} \right)} \right\}\]
\noindent is regarded as the Banach space between ${{L^p}\left( {{\mathbb R^n}} \right)}$ and ${W^{1,p}}\left( {{\mathbb R^n}} \right)$.

\begin{definition}\cite{Podlubny1999FractionalDE} The left sided Riemann-Liouvile fractional integral operator $I_{t}^{\beta } $ of order $\beta >0$, of a function $u\in C_{\kappa } $, $\kappa \ge -1$ is defined as
  \begin{gather*}
      I_{t}^{\beta } u\left(t\right)=\left\{\begin{array}{l} {\frac{1}{\Gamma \left(\beta \right)} \int _{a}^{t}\left(t-\tau \right)^{\beta -1} u\left(\tau \right)d\tau  ,\quad \beta >0,t>0,} \\ {u\left(t\right),{\rm \; \; \; \; \; \; \; \; \; \; \; \; \; \; \; \; \; \; \; \; \; \; \; \; \; \; \; \; \; \;}\quad \; \beta =0.} \end{array}\right.
  \end{gather*}
where  $\Gamma \left( \beta  \right) = \int_0^\infty  {{e^{ - t}}{t^{\beta  - 1}}dt}$ is a Gamma function. 
\end{definition}

\begin{definition}\cite{kilbas2006theory}
The left sided Riemann-Liouvile fractional derivative operator ${\mathbb D}_{t}^{\beta } $ of order $0<\beta <1$, of a function $u\in C$, is defined as
  \begin{gather*}
    {\mathbb D}_{t}^{\beta } u\left(t\right)=\frac{1}{\Gamma \left(1-\beta \right)} \frac{d}{dt} \int _{a}^{t}\left(t-\tau \right)^{-\beta } u\left(\tau \right)d\tau  .  
  \end{gather*}
\end{definition}

\begin{definition}\cite{Podlubny1999FractionalDE} The Caputo time-fractional derivative operator $D_{t}^{\beta } $of order $\beta >0$, is defined as
  \begin{gather*}
     D_{t}^{\beta } u\left(x,t\right)=\frac{\partial ^{n} u\left(x,t\right)}{\partial t^{n} } =\left\{\begin{array}{l} {\frac{1}{\Gamma \left(\beta \right)} \int _{a}^{t}\left(t-\tau \right)^{n-\beta -1} \frac{\partial ^{n} u\left(x,\tau \right)}{\partial \tau ^{n} } d\tau  ,\quad n-1<\beta <n,t>0,} \\ {\frac{\partial ^{n} u\left(t\right)}{\partial t^{n} } ,{\rm \; \; \; \; \; \; \; \; \; \; \; \; \; \; \; \; \; \; \; \; \; \; \; \; \; \; \; \; }\quad \beta =n\in N.} \end{array}\right. 
  \end{gather*}
\end{definition}

\begin{definition}\cite{yu2018time} The one- and two-parameter Mittag-Leffler function is defined as
  \begin{gather*}
     \left\{\begin{array}{l} {E_{\beta } \left(z^{\beta } \right)=\sum _{k=0}^{\infty }\frac{z^{\beta k} }{\Gamma \left(\beta k+1\right)} ,{\rm \; \; \; \; \; }\beta >0 } \\ {E_{\beta ,\gamma } \left(z^{\beta } \right)=\sum _{k=0}^{\infty }\frac{z^{\beta k} }{\Gamma \left(\beta k+\gamma \right)} ,{\rm \; \; }\beta ,\gamma >0 } \end{array}\right.  
  \end{gather*}
\end{definition}

\begin{definition}\cite{galeati2020convergence} For given $f\in L^{2} \left({\mathbb T}^{d} ;{\mathbb R}^{d} \right)$, $f$ is divergence free in sense of if
  \begin{gather*}
      \left\langle f,\nabla g\right\rangle =0{\rm \; \; \; \; }\forall g\in C^{\infty } \left({\mathbb T}^{d} \right).
  \end{gather*}
\end{definition}

\begin{definition}\cite{galeati2020convergence} The orthogonal projection $\Pi $ is given by
   \begin{equation}\label{GrindEQ__2_2_}
    \Pi :f=\sum _{k\in {\mathbb Z}^{d} }f_{k} e_{k}  \mapsto \Pi \, f=\sum _{k\in {\mathbb Z}^{d} }P_{k} f_{k} e_{k}  ,
   \end{equation}
where $P_{k} \in {\mathbb R}^{d} \times {\mathbb R}^{d} $ is the d-dimensional projection on $k^{\bot } $, $P_{k} =I-\frac{k}{\left|k\right|} \otimes \frac{k}{\left|k\right|} \left(k\ne 0\right)$. Moreover, $\Pi _{N} $ is given by
\begin{equation}\label{GrindEQ__2_3_}
f=\sum _{k\in {\mathbb Z}^{d} }f_{k} e_{k}  \mapsto \Pi _{N} \, f=\sum _{k\in {\mathbb Z}^{d} }f_{k} e_{k}  ,
\end{equation}
where $\Pi _{N} :C^{\infty } \left({\mathbb T}^{d} \right)^{{'} } \to C^{\infty } \left({\mathbb T}^{d} \right)$.
\end{definition}

\begin{definition}\cite{niu2012impacts} We say that a random time $\tau$ is a blow-up time (or explosion time ) of the solution $u\left( {t,x} \right)$ to \eqref{GrindEQ__1_1_}, if the following  two conditions are fulfilled:

(\romannumeral1) For any $t < \tau$, $\sup _{x \in \Omega }\left| {u\left( {t,x} \right)} \right| < \infty$ a.s.;

(\romannumeral2) If $\tau  < \infty $, then $\mathop {\lim }\limits_{t \to \tau } {\sup _{x \in \Omega }}\left| {u\left( {t,x} \right)} \right| = \infty$.
\end{definition}

\subsection{Lemma and Property}
We recall here some useful lemma and properties which have been used throughout the paper.
\begin{lemma}\cite{ju2004existence}\label{lemma2.3} Suppose that $q>1$, $p\in \left[q,+\infty \right)$ and $\frac{1}{p} +\frac{\sigma }{2} =\frac{1}{q}$. Suppose that $\Lambda ^{\sigma } f\in L^{q} $, then $f\in L^{p} $ and there is a constant $C\ge 0$ such that
   \begin{gather*}
      \left\| f\right\| _{L^{p} } \le C\left\| \Lambda ^{\sigma } f\right\| _{L^{q} }.  
   \end{gather*}
\end{lemma}
\begin{lemma}\cite{zhou2018weakness}\label{lemma2.4}(Fractional comparison principle) Let $u\left(0\right)=v\left(0\right)$, $u\left(t\right)$ and $v\left(t\right)$satisfies
   \begin{gather*}
       D_{t}^{\beta } u\left(t\right)\ge D_{t}^{\beta } v\left(t\right)
   \end{gather*}
\end{lemma}

\begin{lemma}\cite{zhou2018weakness}\label{lemma2.5}For $u\left(t\right)\ge 0$,
\begin{equation}\label{GrindEQ__2_7_}
D_{t}^{\beta } u\left(t\right)+c_{1} u\left(t\right)\le c_{2} \left(t\right)
\end{equation}
for almost all $t\in \left[0,T\right],$ where $c_{1} >0$ and the function $c_{2} \left(t\right)$ is non-negative and integrable for $t\in \left[0,T\right].$ Then
\begin{equation}\label{GrindEQ__2_8_}
u\left(t\right)\le u\left(0\right)+\frac{1}{\Gamma \left(\beta \right)} \int _{0}^{t}\left(t-s\right)^{\beta -1} c_{2} \left(s\right)ds .
\end{equation}
\end{lemma}

\begin{lemma}\cite{alsaedi2017survey}\label{lemma2.6} For any function $v\left( t \right)$ absolutely continuous on $\left[ {0,T} \right]$, one has the inequality
   \begin{equation}\label{GrindEQ__2_9_}
    v\left( t \right)D_t^\alpha v\left( t \right) \ge \frac{1}{2}D_t^\alpha {v^2}\left( t \right),\quad \quad 0 < \alpha  < 1.
  \end{equation}
\end{lemma}

\begin{lemma}\cite{di2012hitchhikers}\label{lemma2.7} Let $s\in \left(0,1\right)$ and $p\in \left[1,+\infty \right)$ be such that $sp<n.$ Let $\Omega \subseteq {\mathbb R}^{n} $ be an extension domain for $W^{s,p} $. Then there exists a positive constant $C=C\left(n,p,s,\Omega \right)$ such that, for any $f\in W^{s,p} \left(\Omega \right)$, we have
  \begin{equation}\label{GrindEQ__2_10_}
      \left\| u\right\| _{L^{q} \left(\Omega \right)} \le C\left\| u\right\| _{W^{s,p} \left(\Omega \right)} ,
   \end{equation}
for any $q\in \left[p,p^{*} \right]$, where $p^{*} =p^{*} \left(N,s\right)=\frac{Np}{N-sp} $ is the so-called fractional critical exponent; i.e., the space $W^{s,p} \left(\Omega \right)$ is continuously embedded in $L^{q} \left(\Omega \right)$ for any $q\in \left[p,p^{*} \right].$ If, in addition, $\Omega $ is bounded, then the space $W^{s,p} \left(\Omega \right)$ is continuously embedded in $L^{q} \left(\Omega \right)$ for any $q\in \left[1,p^{*} \right].$
\end{lemma}

\begin{lemma}\cite{zhu2018new}\label{lemma2.11} Let $a\left(t\right)$ be a continuous on $\left[0,T\right)\left(0<T\le \infty \right)$, $l\left(t\right)$ is nonnegative and locally integrable on $\left[0,T\right)$, and suppose $u\left(t\right)$ be a continuous nonnegative function on $\left[0,T\right)$ with
   \begin{equation}\label{GrindEQ__2_15_}
         u\left(t\right)\le a\left(t\right)+\int _{0}^{t}l(s)u(s)ds ,t\in \left[0,T\right).
   \end{equation}
Then
    \begin{equation} \label{GrindEQ__2_16_} 
         u\left(t\right)\le a\left(t\right)+\int _{0}^{t}l(s)a\left(s\right)\exp \left(\int _{s}^{t}l\left(\tau \right)d\tau  \right)ds ,t\in \left[0,T\right). 
    \end{equation} 
If $a\left(t\right)$ is a nonnegative non-decreasing on $\left[0,T\right)$, the inequality \eqref{GrindEQ__2_16_} is reduced to
   \begin{equation}\label{GrindEQ__2_17_}
       u\left(t\right)\le a\left(t\right)\exp \left(\int _{0}^{t}l\left(s\right)ds \right)
   \end{equation}
If $a\left(t\right)\equiv 0$, then we can get $u\left(t\right)\equiv 0$ on $\left[0,T\right)$.
\end{lemma}

\begin{lemma}\cite{cortissoz2007skorokhod}\label{lemma2.13}(Skorokhod representation theorem) Suppose\textbf{ }$P_{n} ,\; \; n=1,2,\cdots $ and $P$ are probability measures on $S$ (endowed with its Borel$\sigma $-algebra) such that $P_{n} \Rightarrow P$. Then there is a probability space $\left(\Omega ,{\mathcal F},P\right)$ on which are defined $S$-valued random variables $X_{n} ,\; \; n=1,2,\cdots $ and $X$ with distributions $P_{n} $ and $P$, respectively, such that $\lim _{n\to \infty } X_{n} =X\; \, a.s.$.
\end{lemma}

\begin{lemma}\cite{ren2008burkholder}\label{lemma2.10} For $1\le p<\infty $, there exist constants such that the following holds: for every $N\in {\mathbb N}$ and every martingale $\left(X_{k} \right)_{k=0}^{N} $, we have
   \begin{equation}\label{GrindEQ__2_14_}
        E\left[X\right]_{N}^{{p\mathord{\left/ {\vphantom {p 2}} \right. \kern-\nulldelimiterspace} 2} } \le a_{p} E\left[\left(X_{N}^{*} \right)^{p} \right],\quad \quad E\left[\left(X_{N}^{*} \right)^{p} \right]\le b_{p} E\left[X\right]_{N}^{{p\mathord{\left/ {\vphantom {p 2}} \right. \kern-\nulldelimiterspace} 2} } ,
    \end{equation}
where $X_{N}^{*} :=\mathop{\max }\limits_{k\le N} \left|x_{k} \right|.$
\end{lemma}

\begin{property}\cite{borikhanov2022qualitative}\label{property2.1} If $0<\alpha <1,$$u\in AC^{1} \left[0,T\right]$ or $u\in C^{1} \left[0,T\right],$ then the equality 
\[I_{t}^{\beta } \left(D_{t}^{\beta } u\right)\left(t\right)=u\left(t\right)-u\left(0\right),\] 
and
\[D_{t}^{\beta } \left(I_{t}^{\beta } u\right)\left(t\right)=u\left(t\right),\] 
hold almost everywhere on $\left[0,T\right].$ In addition
\begin{equation}\label{GrindEQ__2_19_}
D_{t}^{1-\beta } \int _{0}^{t}D_{t}^{\beta }  u\left(\tau \right)d\tau =\left(I_{t}^{\beta } \frac{d}{dt} I_{t}^{\beta } I_{t}^{1-\beta } \frac{d}{dt} u\right)\left(t\right)=u\left(t\right)-u\left(0\right).
\end{equation}
\end{property}

\section{Applied literature and motivation} \label{appendixB}

We gather here some updated information on the occurrence of the nonlinear fractional stochastic differential equation we propose and related models in the physical or probabilistic literature.

\begin{enumerate}
    \item[$\bullet$] Stochastic partial differential equations play an increasingly important role in modeling phenomena such as fluid mechanics, atmospheric and oceanography, statistical mechanics, and even finance, cf.\cite{RN85,tien2013fractional,glockle1995fractional}. Existence, uniqueness, regularity, invariant measures and other properties of stochastic partial differential equations have also been widely studied, see \cite{zou2018existence,alsaedi2021global,ju2004existence,RN83}. In recent years, significant advances have also been made in the theory and analysis of fractional-order partial differential equations, which are often used to model anomalous physical processes because of their singularity and non-local character \cite{li2021blow}, to preserve the genetic and memory properties of functions in practical problems \cite{RN89,RN90}. And the physical meaning of the parameters in the fractional models is much clearer.
    \item[$\bullet$] The problem of blow-up or nonexistence of solutions has attracted the interest of many scholars as they are very useful for applications. For deterministic partial differential equations, there exists a vast literature on finite-time explosion, see \cite{galaktionov2002problem} and references therein. In contrast, much less works have been done on such problems for stochastic partial differential equations. Bao and Yuan \cite{bao2016blow} studied the problem of L\'{e}vy-type noise-induced explosions in a class of stochastic reaction-diffusion equations with jumps, Mueller and Sowers \cite{mueller1993blowup} showed critical parameters for blow-up in a class of stochastic heat equations, and subsequently, Muller in \cite{mueller2000critical} further investigated the probability of blow-up in finite time.
    \item[$\bullet$] A wide variety of stochastic partial differential equations that are driven by noise are currently being studied in depth both theoretically and experimentally, because of their convenience in explaining phenomena in many fields, see \cite{GROISMAN2012150,chow2009nonlinear} and references therein. In reality, the noise of a system is diverse and independent, e.g., it is classified into many types, such as deterministic \cite{iyer2021convection} or random, additive or multiplicative, etc., as detailed in \cite{iyer2021convection,GROISMAN2012150,niu2012impacts,chow2011explosive}, In fact, the effect of noise on a system is not always negative. Sometimes, noise may also have a positive effect on bad systems. For example, noise can make an unstable system become a stable stochastic system, see \cite{caraballo2001stabilization}. Mijena and Nane \cite{MIJENA20153301} demonstrated the existence and uniqueness of solutions to equations with temporal white noise. Dozzi \cite{DOZZI2010767} showed that the global solution of the parabolic equation perturbed by temporal white noise is smoother in the spatial variables compared to the equation perturbed by spatially white noise. Inspired by the above works, we discuss whether noise can have a good effect on the explosion time.
\end{enumerate}

\end{appendices}

\bibliography{mybibfile}

\end{document}